\documentclass[12pt]{amsart}

\usepackage{amssymb}
\usepackage{amsthm}
\usepackage{mathptmx}
\usepackage{enumerate}
\usepackage{tikz}
\usepackage{graphicx}
\usepackage[font=footnotesize]{caption}

\usepackage[bookmarks=true]{hyperref}
\hypersetup{colorlinks,citecolor=blue,linkcolor=blue}

\newtheorem{theorem}{Theorem}[section]
\newtheorem{lemma}[theorem]{Lemma}
\newtheorem{prop}[theorem]{Proposition}

\newtheorem{conjecture}[theorem]{Conjecture}

\newtheorem{problem}[theorem]{Problem}
\newtheorem*{quest}{Question}
\theoremstyle{definition}

\newtheorem{example}[theorem]{Example}

\DeclareMathOperator{\SL}{SL}

\DeclareMathOperator{\Ort}{O}

\DeclareMathOperator{\PO}{PO}
\DeclareMathOperator{\SU}{SU}
\DeclareMathOperator{\U}{U}

\DeclareMathOperator{\PSL}{PSL}
\DeclareMathOperator{\PGL}{PGL}

\def\aut{\operatorname{Aut}}
\def\isom{\operatorname{Isom}(\mathbb{H}^n)}

\def\Hy{\mathbb{H}}
\def\Eucln1{\mathbb{E}^{\,n,1}}
\def\Sph{\mathbb{S}}
\def\Mob{\mathrm{Mob}}

\def\Gr{\operatorname{Gr}}

\def\OO{\mathcal{O}}
\def\HG{\mathbb{H}^n/\Gamma}

\def\e_HW{\epsilon_{HW}}

\def\Z{{\mathbb Z}}
\def\Q{{\mathbb Q}}
\def\R{{\mathbb R}}
\def\C{{\mathbb C}}

\def\S{{\mathbb S}}

\def\vol{{\rm Vol}}
\def\cO{{\mathfrak o}_k}

\def\D{\mathcal{D}}

\def\hP{{\mathrm P}}
\def\hF{{\mathrm F}}

\def\G{{\mathrm G}}

\def\Dn{\mathrm{D}}

\begin{document}
\title{Arithmetic hyperbolic reflection groups}

\dedicatory{Dedicated to Ernest Borisovich Vinberg}


\author{Mikhail Belolipetsky}
\address{
IMPA\\
Estrada Dona Castorina, 110\\
22460-320 Rio de Janeiro, Brazil}
\email{mbel@impa.br}

\begin{abstract}
A hyperbolic reflection group is a discrete group generated by reflections in the faces of an $n$-dimensional hyperbolic polyhedron. This survey article is dedicated to the study of arithmetic hyperbolic reflection groups with an emphasis on the results that were obtained in the last ten years and on the open problems.
\end{abstract}

\maketitle

\tableofcontents

\section{Introduction}

Consider a finite volume polyhedron $\hP$ in the $n$-dimensional hyperbolic space $\Hy^n$. It may occur that if we act on $\hP$ by the hyperbolic reflections in its sides, the images would cover the whole space $\Hy^n$ and would not overlap with each other. In this case we say that the transformations form a \emph{hyperbolic reflection group} $\Gamma$ and that $\hP$ is its fundamental polyhedron, also known as the \emph{Coxeter polyhedron} of $\Gamma$. We can give an analogues definition of the spherical and euclidean reflections groups, the classes which are well understood after the work of Coxeter \cite{Cox34}. What kind of properties characterize hyperbolic Coxeter polyhedra? For example, we have to assume that all the dihedral angles of $\hP$ are of the form $\frac{\pi}m$ for some $m \in \{2, 3, \ldots , \infty\}$ because otherwise some images of $\hP$ over $\Gamma$ would overlap. It appears that under some natural assumptions one can prove a general finiteness theorem for the possible types of hyperbolic Coxeter polyhedra in all dimensions. In this paper we shall give an overview of the related methods, results and open problems.

Reflection groups are ubiquitous in mathematics: they appear in group theory, Riemannian geometry, number theory, algebraic geometry, representation theory, singularity theory, low-dimensional topology, and other fields. A vivid description of the history of spherical and euclidean reflection groups can be found in the \emph{Note Historique} of Bourbaki's volume \cite{Bourb_Lie456}. The study of the spherical ones goes back to the mid XIXth century geometrical investigations of M\"obious and Schl\"afli. It then continued in the work of Killing, Cartan and Weyl on Lie theory. In a remarkable paper published in 1934 \cite{Cox34}, Coxeter gave a complete classification of irreducible spherical and euclidean reflection groups. Hyperbolic reflection groups in dimension two were described by Poincar\'e and Dyck already in 1880s \cite{Poin1882, Dyck1882}, they then played a prominent role in the work of Klein and Poincar\'e on discrete groups of isometries of the hyperbolic plane. Much later, in 1970, Andreev proved an analogous result for the hyperbolic three-space giving a classification of convex finite volume polyhedra in $\Hy^3$ \cite{Andreev70a, Andreev70b}. Later on Andreev's theorem played a fundamental role in Thurston's work on geometrization of three-dimensional manifolds. The history and results about reflection groups in algebraic geometry are thoroughly discussed in Dolgachev's survey paper \cite{Dolgachev08_surv} and his lecture notes \cite{Dolgachev15_lect}, the latter giving more details and more focused on the hyperbolic groups. Concluding this very brief overview let us mention that many connections between reflection groups and group theory, combinatorics, and geometry can be found in the book by Conway and Sloane \cite{Conway-Sloane}.

\begin{figure}
\centering
\begin{minipage}{.5\textwidth}
  \centering
  \includegraphics[width=.92\textwidth]{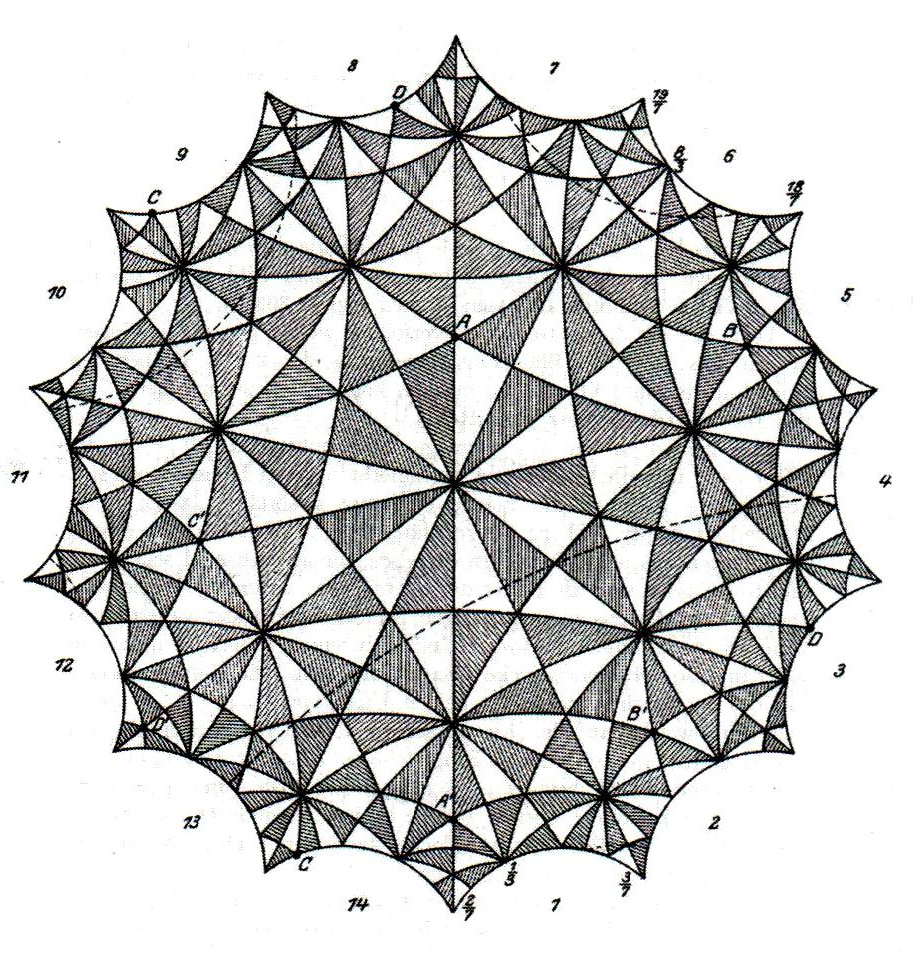}
  \caption{A tiling of a domain on the hyperbolic plane by the $(2,3,7)$-triangles (image taken from Klein's original article \cite{Klein79}).} 
  \label{fig1}
\end{minipage}%
\begin{minipage}{.5\textwidth}
  \centering
  \includegraphics[width=.9\textwidth]{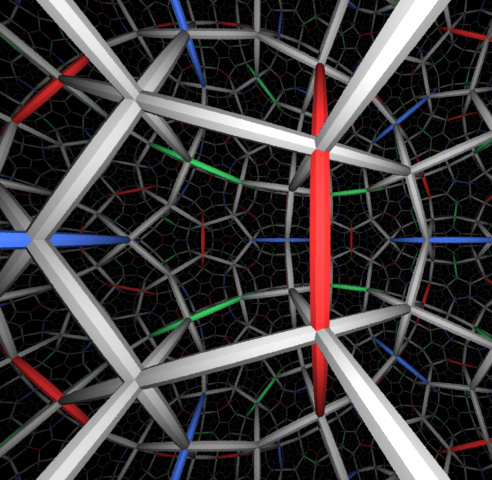}
  \caption{A tiling of the hyperbolic $3$-space by the right-angled dodecahedra (image courtesy of Geometry Center, University of Minnesota).} 
  \label{fig2}
\end{minipage}
\end{figure}

Let us recall some well known examples of hyperbolic reflection groups. Let the dimension $n = 2$,  and consider geodesic triangles $\hP_1$ and $\hP_2$ in the hyperbolic plane with the angles $\frac{\pi}2$, $\frac{\pi}3$, $\frac{\pi}7$ and $\frac{\pi}2$, $\frac{\pi}3$, $0 =\frac{\pi}\infty$, respectively. The corresponding reflection groups $\Gamma_1$ and $\Gamma_2$ are discrete subgroups of the group of isometries of the hyperbolic plane. The first of them is known as the \emph{Hurwitz triangle group}. It is ultimately related with the \emph{Klein quartic surface} $\mathcal{X}$, the corresponding tiling of the fundamental domain of $\mathcal{X}$ on the hyperbolic plane which is shown in Figure~\ref{fig1} appeared in Klein's 1879 paper \cite{Klein79}. This group has many remarkable properties, a fascinating discussion of the related topics can be found in a book \cite{Levy_book1999}. The second group has unbounded fundamental polyhedron whose hyperbolic area is finite (and $=\frac{\pi}6$), it is isomorphic to the \emph{extended modular group} $\PGL(2, \Z)$. Changing to the dimension $n = 3$, we can encounter the right-angled dodecahedron whose corresponding tiling of the hyperbolic three-space as seen from within is represented in Figure~\ref{fig2}. This image was produced by the Geometry Center at the University of Minnesota in the late 1990’s, among other places it appeared in the video ``Not Knot'' available at the Geometry Center homepage \cite{GC} and on the cover of the published edition of the celebrated Thurston's lecture notes \cite{Thur97}. The theories of Klein--Poincar\'e and Thurston were developed for studying much more general classes of discrete groups of isometries, but in both cases hyperbolic reflection groups provided a source of important motivating examples. 

Similar to the way that the polygon in Figure~\ref{fig1} is tiled by triangles, the higher dimensional hyperbolic polyhedra may admit decompositions into smaller parts. A nice computer visualization of some three-dimensional tilings of this kind including two different decompositions of the right-angled dodecahedron is presented in \cite{ACMR09}. In dimensions higher than three it is much harder to draw polyhedra and tilings but we can still study them using different methods. Moreover, the combinatorial and geometrical structure of the fundamental polyhedron of a hyperbolic reflection group in dimension $n\ge 3$ is uniquely determined by its Coxeter diagram (see Section~\ref{sec:prelim}). Throughout this paper we shall study the groups acting in spaces of arbitrary dimensions. 

In what follows we are going to restrict our attention to \emph{arithmetic hyperbolic reflection groups}. The definition of arithmeticity will be given in the next section, here we just note that all the examples considered above are arithmetic. Some results about the general class of hyperbolic reflection groups were discussed in the survey paper by Vinberg \cite{Vinb85}. For the methods that will feature in the present survey arithmeticity plays an essential role.

We now recall a fundamental theorem of Vinberg \cite{Vinb81, Vinb84}, which is well known and at the same time remains surprising:

\begin{theorem}\label{thm1.1}
 There are no arithmetic hyperbolic reflection groups in dimensions $n \ge 30$.
\end{theorem}

One may suspect that the yet to be defined notion of arithmeticity is crucial here but the conjecture is that it is not the case:

\begin{conjecture}
 Theorem \ref{thm1.1} is true without the arithmeticity assumption.
\end{conjecture}

In fact, in his paper Vinberg proved two different theorems one of which says that there are no cocompact hyperbolic reflection groups in dimensions $n \ge 30$, and the other that there are no non-cocompact arithmetic hyperbolic reflection groups in these dimensions. The coincidence of the dimensions in the two theorems is accidental. Later on Prokhorov proved that there are no non-cocompact reflection groups in dimensions $n \ge 996$ thus confirming the conjecture for sufficiently large $n$ \cite{Prokh86}. There is no reason to expect any of these bounds for the dimension to be sharp.
In \cite{Esselmann96}, Esselmann was able to show that the dimensions of non-cocompact arithmetic reflection groups satisfy $n \le 21$, $n \neq 20$. This result is sharp because of the known examples due to Vinberg and Borcherds (see Section~\ref{sec:examples}). Non-arithmetic hyperbolic reflection groups in dimensions bigger than two were first constructed by Makarov and Vinberg in 1960s. The
examples of such groups are currently known for dimensions $n \le 12$, $n = 14$, $18$ \cite{Vinb14}.

In connection with Theorem~\ref{thm1.1} let us mention the following more general question communicated to me by Anton Petrunin:

\begin{quest}
 Do there exist hyperbolic lattices (i.e. discrete cofinite isometry groups) generated by elements of finite order in spaces of large dimension?
\end{quest}
The expected answer to this question is ``No'' but it is far from settled. A related discussion can be found on the MathOverflow page \cite{MOver}. There is a hope that the methods considered in this survey may be applied for attacking this problem for arithmetic lattices, but up to now our attempts to do it were not successful. We shall come back to this question in Section~\ref{sec:problems} dedicated to the open problems.

Another principal question is what happens in the dimensions for which there do exist arithmetic reflection groups --- how to construct the examples of such groups and is it possible to classify all of them? These questions will be in the focus of the discussion in our survey.

\medskip

The content of the paper is as follows. In Section~\ref{sec:prelim} we recall the definition of arithmetic hyperbolic reflection groups and some well known results about them. Section~\ref{sec:Nikulin} contains a brief discussion of the work of Nikulin on finiteness results for hyperbolic reflection groups. In the next section we introduce what we call the spectral method and discuss the finiteness theorems obtained by this method. Some effective results that are obtained by the spectral method are discussed in Sectin~\ref{sec:effective}. The following Section~\ref{sec:classification} is dedicated to what is currently known about classification of arithmetic hyperbolic reflection groups. Section~\ref{sec:examples} presents a collection of examples with an emphasis on those that were discovered after Vinberg's 1985 survey paper was published. Some particularly interesting examples were obtained by Borcherds using modular forms and we dedicate the next section to a discussion of the reflective modular forms. Section~\ref{sec:quasirefl} we consider the reflective quotients, in particular, the so called quasi-reflective groups and also non-reflective groups. Finally, in the last section we discuss open problems.

\medskip

Before starting the paper, let us cite three important articles that provide an overview of the subject from different perspectives. These are the papers by Nikulin \cite{Nikulin81}, Vinberg \cite{Vinb85}, and Dolgachev \cite{Dolgachev08_surv}. The related results were also presented in the ICM talks by Vinberg \cite{Vinberg84_icm} and Nikulin \cite{Nikulin87_icm}. When it is possible, we shall try to minimize the overlap with these papers and to focus our attention on the results that were obtained in the last 10 years.

\section{Quadratic forms and arithmetic reflection groups} \label{sec:prelim}

\def\ia{{\mathfrak a}}

Consider an $(n+1)$-dimensional vector space $\Eucln1$ with the inner product defined by a quadratic form $f$ of signature $(n, 1)$. Let
\[\{v \in \Eucln1 | (v, v) < 0\} = \mathfrak{C} \cup (-\mathfrak{C}),\]
where $\mathfrak{C}$ is an open convex cone. In the \emph{vector model}, the hyperbolic space $\Hy^n$ is identified with the set of rays through the origin in $\mathfrak{C}$, or $\mathfrak{C}/\mathbb{R}^+$, so that the isometries of $\Hy^n$ are the orthogonal transformations of $\Eucln1$. We refer to \cite{Geometry2_EMS_part1} and \cite{Ratcliffe_book} for a detailed study of the properties of the vector model.

Let $k$ be a totally real number field with the ring of integers $\cO$, and let $f$ be a quadratic form of signature $(n,1)$ defined over $k$ and such that for every non-identity embedding $\sigma :k\to\R$ the form $f^\sigma$ is positive definite. The group $\Gamma = \Ort_0(f,\cO)$ of the integral automorphisms of $f$ is a discrete subgroup of $H = \Ort_0(n,1)$, which is the full group of isometries of the hyperbolic $n$-space $\Hy^n$ (the group $\Ort_0(n,1)$ is the subgroup of the orthogonal group $\Ort(n,1)$ that preserves the cone $\mathfrak{C}$). The discreteness of $\Gamma$ follows from the discreteness of $\cO$ in $\R^{[k:\Q]}$ and compactness of $\Ort_0(f^\sigma)$ for $\sigma \neq \sigma_{id}$. Using reduction theory, it can be shown that the covolume of $\Gamma$ in $H$ is finite --- this is a special case of the fundamental theorem of Borel and Harish-Chandra \cite{BorHC62} (we refer to the book \cite{OnVinb00} for an accessible exposition of the main ideas of the proof). Discrete subgroups of finite covolume are called \emph{lattices}. The groups $\Gamma$ obtained in this way and subgroups of $H$ which are commensurable with them are called \emph{arithmetic lattices of the simplest type}. The field $k$ is called the \emph{field of definition} of $\Gamma$ (and subgroups commensurable with it). We shall also apply this terminology to the corresponding quotient orbifolds $\Hy^n/\Gamma$.

There are compact and finite volume non-compact arithmetic quotient spaces. A theorem known as the \emph{Godement's compactness criterion} implies that $\Gamma$ is non-cocompact if and only if $\Ort_0(f,k)$ has a non-trivial unipotent element \cite{BorHC62} (see also \cite[Chapter~3, Section~3.3]{OnVinb00}). This is equivalent to the condition that $k = \Q$ and the form $f$ is isotropic. The Hasse--Minkowski theorem implies that for $k = \Q$ and $n \ge 4$ the latter condition automatically holds. Therefore, for $n \ge 4$ the quotient $\Hy^n/\Gamma$ is non-compact if and only if $\Gamma$ is defined over the rationals. For $n = 2$ and $3$ the non-cocompact subgroups are still defined over $\Q$ but there also exist cocompact arithmetic subgroups with the same field of definition.

In general arithmetic subgroups of semisimple Lie groups are defined using algebraic groups. This way the classification of semisimple algebraic groups \cite{Tits66} implies a classification of the possible types of arithmetic subgroups. It follows from the classification that for hyperbolic spaces of even dimension all arithmetically defined subgroups are arithmetic subgroups of the simplest type. For odd $n$ there is another family of arithmetic subgroups given as the groups of units of appropriate Hermitian forms over quaternion algebras. Moreover, if $n=7$ there is also the third type of arithmetic subgroups of $H$ which are associated to the Cayley algebra. The following lemma of Vinberg shows that for our purpose it will be always sufficient to consider only arithmetic subgroups of the simplest type:

\begin{lemma}\cite[Lemma~7]{Vinb67}\label{lem:vinb67}
Any arithmetic lattice $\Gamma\subset H$ generated by reflections is an arithmetic lattice of the simplest type.
\end{lemma}

A discrete subgroup of a Lie group $H$ is called \emph{maximal} if it is not properly contained in any other discrete subgroup of $H$. It is well known that in semisimple Lie groups any lattice is contained in a maximal lattice. Let $\Gamma$ be an arithmetic subgroup which is commensurable with $\G(\cO)=\Ort(f,\cO)$ for some quadratic form $f$ defined over $k$. There exists an $\cO$-lattice $L$ in $k^{n+1}$ such that $\Gamma\cap \G(k) \subset \G^L = \{g\in \G(k) \mid g(L) = L\}$ (cf. \cite{Vinb71}). If $\cO$ is a principal ideal domain in an appropriate basis the transformations from $\G^L$ can be written down by matrices with elements in $\cO$. By Theorem~5 from \cite{Vinb71}, if $\Lambda < \Gamma$ is an arithmetic subgroup generated by reflections then $\Lambda$ is definable over $\cO$. We refer to \cite{Vinb71} for more information about the rings of definition of arithmetic subgroups. Note that the term ``field of definition'' is used by Vinberg in a more restricted sense.


A certain subclass of arithmetic subgroups will play a special role in our study. An arithmetic subgroup $\Gamma < \Ort_0(f)$ is called a \emph{congruence subgroup} if there exists a nonzero ideal $\ia\subset\cO$ such that $\Gamma\supset \Ort_0(f,\ia)$, where $\Ort_0(f,\ia) = \{ g\in\Ort_0(f,\cO) \mid g\equiv \mathrm{Id\;(mod\;\ia)}\}$ is the \emph{principal congruence} subgroup of $\Ort_0(f,\cO)$ of level $\ia$. Congruence subgroups possess a remarkable set of special geometric and algebraic properties. In particular, the first non-zero eigenvalue of the Laplacian on $\Hy^n/\Gamma$ is bounded away from zero by a constant which depends only on $n$ but not on $\Gamma$. This \emph{spectral gap} property will play a key role for the methods that we are going to discuss in Sections~\ref{sec:spec_method}--\ref{sec:effective}.

We have the following fact connecting maximal and congruence arithmetic subgroups:
\begin{lemma}\label{lem:max_ar}
Maximal arithmetic subgroups are congruence.
\end{lemma}
For the arithmetic subgroups defined by quadratic forms the proof of the lemma can be found in \cite[Lemma~4.7]{ABSW08}, where it is based on a material from \cite{PlaRap94}.

\medskip

An interested reader can find much more information about arithmetic groups and their properties in the books \cite{PlaRap94} and \cite{WitMor15}.

\medskip

Let us now come back to the reflection groups --- the objects of our primary interest. A reflection group $\Gamma$ is called a \emph{maximal reflection group} if there is no other discrete subgroup $\Gamma' < \isom$ such that $\Gamma < \Gamma'$ and $\Gamma'$ is generated by reflections. Maximal reflection groups are not necessarily maximal lattices but there is a relation between the two, and it is captured by the following lemma due to Vinberg:
\begin{lemma}\label{lem:max_refl}
A maximal reflection group $\Gamma$ is a normal subgroup of a maximal lattice $\Gamma_{0}$. Moreover, there is a finite subgroup $\Theta<\Gamma_{0}$
such that $\Theta \to  \Gamma_{0}/\Gamma$ is an isomorphism, and $\Theta$ is the group of symmetries of the Coxeter polyhedron of $\Gamma$.
\end{lemma}
This lemma was proved in \cite{Vinb67}, the argument is also reproduced in \cite{ABSW08}.

\medskip

Given an admissible quadratic form $f$ as above, we would like to know when the arithmetic subgroup $\Gamma = \Ort_0(f,\cO)$ is, up to finite index, generated by hyperbolic reflections. If this is the case, the form $f$ and the group $\Gamma$ are called \emph{reflective}. The main practical tool for deciding reflectivity is \emph{Vinberg's algorithm} \cite{Vinb72} which we are going to review now.

In the vector model of $\Hy^n$, a hyperplane is given by the set of rays in $\mathfrak{C}$ which are orthogonal to a vector $e$ of positive square in $\Eucln1$. A hyperplane $\Pi_e$ defines two halfspaces, $\Pi_e^+$ and $\Pi_e^-$, where $'\pm'$ is the sign of $(e,x)$ for $x$ in the corresponding halfspace, and a \emph{reflection}
\begin{equation*}
R_e: x \to x - 2\frac{(e,x)}{(e,e)}e,
\end{equation*}
where the inner product $(u,v) = \frac{1}{2}(f(u+v)-f(u)-f(v))$ is induced by $f$.

We shall assume that $\cO$ is a principal ideal domain (PID). The vector $e$ corresponding to the reflection $R_e$ is defined up to scaling, so if $e$ has $k$-rational coordinates we can normalize it so that the coordinates are coprime integers in $\cO$. With this normalization we can assign to $R_e$ a parameter $s = (e,e)\in \cO$ and call $R_e$ an \emph{$s$-reflection}. 
The reflection $R_e$ belongs to the group $\Ort_0(f,\cO)$ if
$\frac{2}{s}(e,v_i) \in \cO, \text{ for the standard basis vectors } v_i,\ i = 0,\ldots,n$.
Following \cite{Vinb72}, we call this the \emph{crystallographic condition}.

We begin by considering the stabilizer subgroup of a vector $u_0$ with integral coordinates which corresponds to a point $x_0\in \overline{\Hy^n}$ (when dealing with non-cocompact lattices it may be convenient to choose an ideal point $x_0\in \partial\Hy^n$ as a starting node). Consider the (finite) group generated by all reflections in $\Gamma$ whose mirrors pass through $x_0$. Let
$$\hP_0 = \bigcap_{i=1}^m \Pi_{e_i}^-$$
be a fundamental chamber of this group. All the halfspaces $\Pi_{e_i}^-$ are \emph{essential} (i.e. not containing the intersection of the other halfspaces). The corresponding vectors $e_i$ satisfy $(e_i, e_i) > 0$, $(e_i, u_0) = 0$ for all $i$, and the reflections $R_{e_i}$ generate the stabilizer of $x_0$ in $\Ort_0(f; \cO)$. There is a unique fundamental polyhedron $\hP$ of the reflection subgroup which sits inside $\hP_0$ and contains $x_0$. The point $x_0$ is not necessarily a vertex of $\hP$ but it is usually convenient to choose $u_0$ in such a way that $x_0$ is a vertex.

The algorithm continues by picking up further $\Pi_{e_i}$ so that
$$\hP \subseteq \bigcap_{i} \Pi_{e_i}^-.$$
This is done by choosing $e_i$ satisfying the crystallographic condition such that $(e_i, e_i) > 0$, $(e_i, u_0)<0$, $(e_i, e_j) \le 0$ for all $j < i$, and the distance between $x_0$ and $\Pi_{e_i}$ is the smallest possible, i.e. minimizing the value
$$\sinh^2\left(\mathrm{dist}(x_0, \Pi_{e_i})\right) = -\frac{(e_i,u_0)^2}{(e_i,e_i)(u_0,u_0)}.$$
The latter condition implies that all the hyperplanes $\Pi_{e_i}$ are essential. Note that if $k \neq \Q$, its integers do not form a discrete subset of $\R$. Bugaenko showed that regardless of this, the arithmeticity assumption implies that the set of distances considered above is discrete and hence we can always choose the smallest one (see \cite{Bug84}, \cite{Bug90}, \cite{Bug92} or \cite{Mark_thesis}).

The algorithm terminates if it generates a configuration $\hP = \bigcap_{i} \Pi_{e_i}^-$ that has finite volume, in which case the form $f$ is reflective. The finite volume condition can be effectively checked from the Coxeter diagram of $\hP$ --- it is equivalent to each edge of the polyhedron having two vertices, either one or both of which may be at the ideal boundary of the hyperbolic space.

Let us recall that the fundamental polyhedra of the reflection groups are usually described using \emph{Coxeter diagrams}. These are the graphs with vertices corresponding to the vectors $e_i$ (equivalently, the faces of $\hP$). Two different vertices $e_i$, $e_j$ are connected by a \emph{thin edge} of integer weight $m_{ij} \ge 3$ or by $m_{ij}-2$ edges if the corresponding faces intersect with the dihedral angle $\frac{\pi}{m_{ij}}$, by a \emph{thick edge} if they intersect at infinity (dihedral angle zero), and by a \emph{dashed edge} if they are divergent. In particular, two vertices are not joined by an edge if and only if the corresponding faces of $\hP$ are orthogonal. Note that there are some differences between the Coxeter diagrams and the Dynkin diagrams which are used in Lie theory. In particular, the triple edge in our labeling convention means the angle $\frac\pi5$, while on the Dynkin diagram of a Weyl chamber it corresponds to the angle $\frac\pi6$.

In order to produce the Coxeter diagram of the polyhedron $\hP$ we can compute the dihedral angles between the intersecting hyperplanes from the standard formula
$$\cos\left(\frac{\pi}{m_{ij}}\right) = \frac{-(e_i, e_j)}{\sqrt{(e_i,e_i)(e_j,e_j)}}.$$

\begin{example} Let us consider the quadratic form
$$f = 2\left(x_1^2+x_2^2 +x_3^2+x_4^2-x_1x_2-\frac12(1+\sqrt5)x_2x_3-x_3x_4\right).$$
It is an admissible even quadratic form of the discriminant $3-2\sqrt5$ defined over the field $k = \Q(\sqrt{5})$.

We can start running algorithm with $u_0 = (\frac32, 3, 2\phi, \phi)$, where $\phi = \frac{1 + \sqrt 5}2$ is the fundamental unit of $k$, so that the stabilizer subgroup of the corresponding point $x_0\in \Hy^3$ is generated by reflections corresponding to the vectors
$$e_1 = (1,0,0,0),\ e_2 = (0,0,1,0)\text{ and } e_3 = (0,0,0,1).$$
The algorithm finds the forth vector $e_4 = (-1,-1,-\phi,-\phi)$ and terminates.

\newpage

The Coxeter diagram of the resulting configuration is
\begin{figure}[!ht]
  \centering
  \begin{tikzpicture}
    \draw[thick] (0,0) -- (1, 0);
    \draw[thick, double distance = 3pt] (1, 0) -- (2, 0);
    \draw[thick] (1,0) -- (2, 0);
    \draw[thick] (2,0) -- (3, 0);
    \filldraw[fill=white] (0, 0) circle (2pt) node [above] {\scriptsize 1};
    \filldraw[fill=white] (1, 0) circle (2pt) node [above] {\scriptsize 4};
    \filldraw[fill=white] (2, 0) circle (2pt) node [above] {\scriptsize 3};
    \filldraw[fill=white] (3, 0) circle (2pt) node [above] {\scriptsize 2};
  \end{tikzpicture}
\end{figure}

The polyhedron $\hP$ is a bounded simplex in $\Hy^3$ and the group $\Gamma$ generated by reflections in its sides is a well known arithmetic lattice. A difficult theorem of Gehring, Martin and Marshall \cite{GehrMart09, MarshMart12} shows that the order $2$ extension of $\Gamma$ is the minimal covolume lattice in $\mathrm{Isom}(\Hy^3)$, hence we can think of it as the three dimensional analogue of the $(2,3,7)$-Hurwitz group.

Note that the form $f$ is not diagonalizable over $\cO$ and it can be verified that the group $\Gamma$ cannot be obtained as a reflection subgroup of the group of units of some diagonal quadratic form.
\end{example}

We refer to \cite{Vinb72, VK78, Bug84, Bug90, shaiheev90, Bug92, SchWal92, Nikulin00, Allc12, BelMcl13, mcleod_thesis, Mark15, Mark_thesis} for many other examples of the Coxeter polyhedra and their diagrams produced by the algorithm. In most of these papers the form $f$ is defined over $\Q$ which implies (for $n \ge 4$) that the resulting Coxeter polyhedra have cusps. The exception are the Bugaenko's papers where the algorithm was first applied to the forms with coefficients in the real quadratic fields leading to examples of cocompact arithmetic hyperbolic reflection groups in dimensions $n \le 8$. We shall come back to the discussion of the known examples in Section~\ref{sec:examples}.

Let us point out that Vinberg's algorithm has an unfortunate property that it \emph{never halts if the form is not reflective}. Fortunately, in practice this problem can be often bypassed: if a computer implementation of the algorithm produces, say, more than $10^3$ generators for the reflection subgroup, we can expect that it is not a lattice. The latter condition can be rigorously checked by the group theoretic methods in each particular case, for example, by detecting an infinite order symmetry of the reflection polyhedron. It then implies that the group was not reflective. Many concrete examples and some general methods for this verification can be found in the above cited papers.

\section{The work of Nikulin} \label{sec:Nikulin}

In a series of papers beginning with the 1980 article \cite{Nikulin80}, Nikulin established finiteness of the number of commensurability classes of arithmetic hyperbolic reflection groups and obtained upper bounds for degrees of their fields of definition. In this section, we shall briefly review Nikulin's method.

Let $\hP$ be a convex polyhedron in $\Hy^n$. It is an intersection of finitely many halfspaces $\Pi_{e_i}^-$, where the vectors $e_i$ orthogonal to the faces of $\hP$ can be chosen to have square $2$ and directed outward. The matrix
$$A(P) = (a_{ij}) = (e_i, e_j)$$
is called the \emph{Gram matrix} of $\hP$. It uniquely determines $\hP$ up to motions of the ambient space $\Hy^n$.

The polyhedron $\hP$ is a fundamental polyhedron of a discrete reflection group in $\Hy^n$ if and only if $a_{ij} \le 0$ and $a_{ij} = -2\cos(\frac\pi{m_{ij}})$ where $m_{ij} \ge 2$ is an integer whenever $a_{ij} > -2$ for all $i \neq j$. Symmetric real matrices $A$ satisfying these conditions and having all their diagonal elements equal to $2$ are called the \emph{fundamental matrices}.

\medskip

\noindent{\bf Note.} In his papers, Nikulin uses the opposite sign convention that is more common in algebraic geometry. Here we keep the notation that was introduced in the previous section, which is also consistent with the one used by Vinberg.

\medskip

Given a real $t > 0$, we say that the fundamental Gram matrix $A = (a_{ij})$ and the corresponding polyhedron $\hP$ \emph{have minimality $t$} if $|a_{ij}| < t$ for all $a_{ij}$. We can analogously define the \emph{minimality of a face} $\hF$ of $\hP$ by considering the Gram matrix $A(\hF)$ formed by the inner products of the vectors associated to the faces of $\hP$ that have non-trivial intersections with $\hF$ in $\overline{\Hy^n}$. The notion of minimality is central for Nikulin's method.

Suppose that $\hP$ is a fundamental polyhedron of an arithmetic reflection group $\Gamma = \Gamma(\hP)$ in $\Hy^n$. Vinberg \cite{Vinb67} proved that for a finite volume polyhedron $\hP$ it is equivalent to the conditions that all the cyclic products
$$b_{i_1\ldots i_m} = a_{i_1i_2}a_{i_2i_3}\cdots a_{i_{m-1}i_m}a_{i_mi_1}$$
are algebraic integers, the field $K = \Q(\{a_{ij}\})$ is totally real, and the matrices $A^\sigma = (a_{ij}^\sigma)$ are nonnegative definite for every embedding $\sigma: K \to \R$ not equal to identity on the field of definition $k = \Q(\{b_{i_1\ldots i_m}\})$, which is generated by the cyclic products.

A fundamental matrix $A(\hP)$ and the corresponding reflection group $\Gamma(\hP)$ is called \emph{$V$\mbox{-}\nobreak\hspace{0pt}arithmetic} if it satisfies the conditions of Vinberg's arithmeticity criterion except that $\hP$ is not required to have finite volume. The property of $V$-arithmeticity is much easier to check than arithmeticity and, moreover, it is inherited by sub-polyhedra: if $\hP'$ is an intersection of a subset of the halfspaces $\Pi_{e_i}^-$ defining the polyhedron $\hP$ and its Gram matrix $A(\hP')$ is indefinite than $\hP'$ is also $V$-arithmetic with the same field of definition (also called \emph{ground field}) as $\hP$.

The hereditary property of $V$-arithmeticity allows to reduce some questions to the analysis of simple configurations. The basic case is called an \emph{edge polyhedron (chamber)}. It refers to a fundamental matrix $A(\hP)$ such that all the corresponding hyperplanes $\Pi_e$ contain at least one of the two distinct vertices $v_1$, $v_2$ of a one-dimensional edge $v_1v_2$ of $\hP$. If the vertices $v_1$, $v_2$ are finite, the edge chamber is called \emph{finite}. In this part we are mainly interested in bounding the degree of the fields of definition of arithmetic reflection groups, thus we can restrict to the finite case. With the necessary modifications, Nikulin's method also applies to the unbounded Coxeter polyhedra. Any finite edge chamber has precisely $n+1$ sides corresponding to the vectors $e_1$, $e_2$ and $n-1$ vectors $\{ e_j \}_{j\in J}$ whose hyperplanes contain the full edge $v_1v_2$. The Gram matrix $A(\hP)$ is hyperbolic but its submatrices corresponding to $e_1 \cup \{ e_j \}_{j\in J}$ and $e_2 \cup \{ e_j \}_{j\in J}$ are positive definite. The only element of $A(\hP)$ that can have absolute value bigger than $2$ is $u = (e_1, e_2)$. It follows that an edge polyhedron has minimality $t > 2$ if and only if $|u| = |(e_1, e_2)| < t$. The Coxeter diagram $G$ of an edge chamber $\hP$ has exactly one connected component $G(\hP^{hyp})$ corresponding to a hyperbolic submatrix of $A(\hP)$ and containing $e_1$ and $e_2$, and possibly several positive definite connected components. The Gram matrix $A(\hP^{hyp})$ corresponds to an edge chamber of dimension $\# G(\hP^{hyp})-1$. If $\hP$ is $V$-arithmetic, then $\hP$ and the hyperbolic connected component $\hP^{hyp}$ are defined over the same field $k$.

We now can state a principal technical theorem of Nikulin:

\begin{theorem}\cite[Theorem~2.3.1]{Nikulin81}\label{thm:nikulin1}
Given any $t > 0$, there is an effective constant $N(t)$ such that for every $V$-arithmetic edge chamber of minimality $t$ with the ground field $k$ of degree more than $N(t)$ over $\Q$, the number of vertices in the hyperbolic connected component of the Coxeter graph is less than $4$.
\end{theorem}

The proof of Theorem~\ref{thm:nikulin1} in \cite{Nikulin81} uses a variant of Fekete's theorem (1923) on the existence of non-zero integer polynomials of bounded degree with small deviation from zero on appropriate intervals (see also \cite[Section 6]{Nikulin11} for a review of the proof and some corrections).

The minimality $t = 14$ is especially important for Nikulin's method. This is because a fundamental polyhedron of an arithmetic hyperbolic reflection group (not assumed to be cocompact here) always has a face with minimality $14$ --- a result that was proved in the early Nikulin's papers (see \cite[Lemma~3.2.1]{Nikulin80} and \cite[proof of Theorem~4.1.1]{Nikulin81}). This fact allows one to reduce various finiteness and classification problems to the estimation of the value of the \emph{transition constant} $N(14)$. This way in \cite{Nikulin07_fin}, Nikulin showed that his previous work together with the finiteness results for dimensions $n = 2$ and $3$ obtained by the spectral method, which we are going to discuss in the next section, can be applied to prove a general finiteness theorem for the number of commensurability classes of arithmetic hyperbolic reflection groups. He also proved that the degrees of the fields of definition of arithmetic hyperbolic reflection groups are bounded above by the maximum of $N(14)$ and the maximal possible degree in dimensions two and three. In the subsequent papers Nikulin gave explicit bounds for the constant $N(14)$, for example, in \cite{Nikulin09_flds1} (see also Remark~5.1 in \cite{Nikulin11} for a correction) he showed that $N(14) \le 120$. The best result of this kind is obtained in \cite{Nikulin11}, where it is shown that $N(14) \le 25$. It is worth mentioning that most of these results make use of Theorem~\ref{thm:nikulin1} and its proof. The main difference between the latest improvement and the previous papers is that there, instead of relying only on estimates for the Fekete's existence theorem, Nikulin constructed certain explicit polynomials with the required properties.

The result of \cite{Nikulin11}, together with \cite{Macl11} and \cite{BelLin14} which we shall discuss later on in Section~\ref{sec:effective}, implies:
\begin{theorem} (cf. \cite[Corollary~5.1]{BelLin14})\label{thm:nikulin2}
The degree of the fields of definition of arithmetic hyperbolic reflection groups in all dimensions is at most $25$.
\end{theorem}
This is currently the best known general bound.

\section{Spectral method and finiteness theorems} \label{sec:spec_method}

Given a lattice $\Gamma \le \isom$, we have an associated quotient Riemannian orbifold $\mathcal{O} = \Hy^n/\Gamma$. If $\Gamma$ is generated by reflections we call $\mathcal{O}$ a \emph{reflection orbifold}. The main idea of the method discussed in this section is that the global geometric properties of the reflection orbifolds can provide us important information about the hyperbolic reflection groups.

The spectral method that we are going to review is based on some properties of the spectrum of the Laplacian on orbifolds. It was first applied for proving finiteness results for arithmetic hyperbolic reflection groups by Long, Maclachlan and Reid \cite{LMR06} and by Agol \cite{Agol06icm} in dimensions $n = 2$ and $3$, respectively. The first paper makes use of the Zograf's spectral proof of Rademacher's conjecture for congruence subgroups of the modular group \cite{Zograf91}. In his work \cite{Agol06icm}, Agol found that one can employ the Li--Yau inequality for conformal volume instead of Zograf's spectral inequality for surfaces. This allows to extend the domain of applicability of the method to a much wider class of spaces. Agol completed the proof of the finiteness theorem for $n = 3$; in a later joint work with Agol, Storm and Whyte we showed how to extend the argument to an arbitrary dimension \cite{ABSW08}. We shall now review this approach.

Conformal volume was introduced by Li and Yau in \cite{LiYau82}, partially motivated by generalizing results on surfaces due to Yang--Yau \cite{YangYau80}, Hersch \cite{Hersch70}, and Szeg{\"o} \cite{Szego54}. In \cite{ABSW08}, we generalized this notion to orbifolds. Let $(\mathcal{O},g)$ be a complete Riemannian orbifold, possibly with boundary and let $|\mathcal{O}|$ denote the underlying topological space. Denote the volume form by $dv_{g}$, and the volume by $\vol(\mathcal{O},g)$. Let $\Mob(\Sph^{n})$ denote the group of conformal transformations of $\Sph^{n}$. It is well known that $\Mob(\Sph^{n})=\mathrm{Isom}(\Hy^{n+1})$. The topological space $|\mathcal{O}|$ has a dense open subset which is a Riemannian manifold. We call a map $\varphi: |\mathcal{O}_{1}|  \to |\mathcal{O}_{2}|$ a {\it PC map} if it is a continuous map which is piecewise a conformal immersion. Clearly, if $\varphi:|\mathcal{O}|\to \Sph^{n}$ is PC, and $\mu\in \Mob(\Sph^{n})$, then $\mu\circ \varphi$ is also a PC map.

Let $(\Sph^{m},can)$ be the $m$-dimensional sphere with the canonical round metric. For a piecewise smooth map $\varphi: |\mathcal{O}| \to (\Sph^{m},can)$, define
$$V_{PC}(m,\varphi) = \underset{\mu\in \Mob(\Sph^{m})}{\sup} \vol(\mathcal{O}, (\mu\circ \varphi)^{*}(can)).$$
If there exists a PC map $\varphi: |\mathcal{O}| \to \Sph^{m}$, then we also define
$$V_{PC}(m, \mathcal{O}) = \underset{\varphi:|\mathcal{O}|\to \Sph^{m}  PC}{\inf} V_{PC}(m,\varphi).$$
We call $V_{PC}(m,\mathcal{O})$ the $m$-dimensional {\it piecewise conformal volume} of $\mathcal{O}$.

Using the Nash embedding theorem, it can be shown that the $m$-dimensional conformal volume is always well defined for a sufficiently large $m$. It is clear that $V_{PC}(m,\mathcal{O})\geq V_{PC}(m+1,\mathcal{O})$, hence we can define
the (piecewise) \emph{conformal volume}
$$\vol_{c}(\mathcal{O})=\underset{n\to \infty}{\lim} V_{PC}(n,\mathcal{O}).$$
We refer to \cite{ABSW08} for further discussion and basic properties of the conformal volume. One of the immediate corollaries of the definitions allows us to compute the conformal volume of the reflective orbifolds (cf. Facts 3, 4 in \cite[Section 2]{ABSW08}):

\begin{quote}\emph{
Let $\Gamma_r < \isom$ be a lattice generated by reflections in hyperplanes in $\Hy^n$ and let $\OO_r = \Hy^n/\Gamma_r$. Then we have
 $
 \vol_c(\OO_r) = \vol(\Sph^n, can).
 $
}\end{quote}
Now the Li--Yau inequality generalized to $n$-dimensional orbifolds says \cite{ABSW08}:
\begin{equation}\label{eq:Li-Yau}
 \lambda_1(\OO) \cdot \vol(\OO)^\frac{2}{n} \le n\cdot \vol_c(\OO)^\frac{2}{n}.
\end{equation}
Here $\lambda_1(\OO)$ denotes the first nonzero eigenvalue of the Laplacian on $\OO$ (called the \emph{spectral gap}), $\vol$ is the hyperbolic volume and $\vol_c$ is the conformal volume defined above.

Equation \eqref{eq:Li-Yau} shows that information about the spectral gap and hyperbolic volume of a reflection orbifold can be played against the upper bound given by the value of its conformal volume. The information that we need here can be deduced from arithmeticity. Indeed, if $\Gamma$ is a congruence subgroup of $\isom$ (cf. Section~\ref{sec:prelim} for the terminology), then the well known conjectures of Ramanujan and Selberg imply that
$$\lambda_1(\HG) \ge \lambda_1(\Hy^n).$$
These conjectures are still open but less precise low bounds are known, and for our purpose they can serve almost as good as the conjectures. We have:
\begin{equation*}
\lambda_1(\HG) \ge \delta(n),
\end{equation*}
where $\delta(2) = \frac{3}{16}$ by Vigneras~\cite{Vign83} and if $n \ge 3$, $\delta(n) = \frac{2n-3}{4}$ by Burger--Sarnak~\cite{BurgSarn91}. Moreover, if $\Gamma$ is defined by a quadratic form (which is always the case for the arithmetic hyperbolic reflection groups by Lemma~\ref{lem:vinb67}), then more recent work of Luo, Rudnick and Sarnak implies that we can take $\delta(2) = 0.21$ and $\delta(n) = \frac{15n-24}{25}$ for $n\ge 3$ \cite{LuoRS99}. The proofs of these bounds are based on deep results about automorphic representations.

Now let us assume that $\Gamma$ is at the same time a congruence subgroup and a reflection group, and let $\OO = \HG$. Following the argument in \cite{Bel11}, we can then quickly prove the two principal finiteness theorems. We have:
\begin{equation*}
\delta(n)\cdot\vol(\OO)^\frac2n \le n\cdot\vol(\Sph^n)^\frac2n;
\end{equation*}
\begin{equation}\label{eq:pr2}
\vol(\OO) \le \left(\frac{n}{\delta(n)}\right)^\frac{n}{2}\vol(\Sph^n).
\end{equation}

By the theorems of Wang~\cite{Wang72} for $n\ge4$ and Borel~\cite{Bor81} for $n = 2,\ 3$ there are only
finitely many (up to conjugacy) arithmetic subgroups of $\isom$ of bounded covolume. As the
right-hand side of \eqref{eq:pr2} depends only on the dimension, we immediately obtain our first
finiteness theorem:

\begin{theorem} \label{thm4.1}
For every $n\ge2$ there are only finitely many conjugacy classes of congruence reflection subgroups of $\isom$.
\end{theorem}

Let $n\ge 3$. We have $\delta(n)\ge\frac{15n-24}{25}$;
\begin{equation}\label{eq:pr3}
\frac{\vol(\OO)}{\vol(\S^n)}\le \left(\frac{25n}{15n-24}\right)^\frac{n}{2} \le \left(\frac{25}{7}\right)^\frac{n}{2}.
\end{equation}

At this place we need to recall some recent results about volumes of arithmetic hyperbolic $n$-orbifolds. By \cite{Bel04}, \cite{Bel07} and \cite{BelEme}, $\vol(\HG)$ is bounded below by a function which grows super-expo\-nen\-tial\-ly with $n$. These results are discussed in some detail in my other survey article \cite{Bel14icm}, the required corollary in a more precise form can be found as Corollary~3.3 there. As $\vol(\Sph^n)\to0$ when $n\to\infty$, the super-exponential lower bound holds true also for the quotient $\vol(\HG)/\vol(\Sph^n)$. Hence the left-hand side of \eqref{eq:pr3} grows super-exponentially with $n$ while the right-hand side is only exponential. This gives our second finiteness theorem:

\begin{theorem}\label{thm4.2}
If $n$ is sufficiently large then $\isom$ does not contain any congruence reflection subgroups.
\end{theorem}

\medskip

The main problem with this argument is that we cannot expect that all maximal arithmetic reflection groups are congruence --- the counterexamples are known in dimensions $2$ and $3$ by the work of Lakeland \cite{LakelandPhD, Lakeland12}, and it is not clear what happens in higher dimensions. The goal of \cite{ABSW08} was to show how finer arithmetic techniques can be applied in order to partially bypass this difficulty and prove an analogue of Theorem~\ref{thm4.1} for all maximal arithmetic reflection groups. We can not yet prove an analogue of Theorem~\ref{thm4.2} using the spectral method, but we know that it is true thanks to the previous work of Vinberg (cf. Theorem~\ref{thm1.1} in the introduction). Thus we have:

\begin{theorem} \label{thm4.3}
There are only finitely many conjugacy classes of arithmetic maximal hyperbolic reflection groups.
\end{theorem}

This result was proved independently in \cite{ABSW08} and \cite{Nikulin07_fin}. It implies that in principle it is possible to give a complete classification of the arithmetic hyperbolic reflection groups.

Let us note that neither arithmeticity nor maximality assumptions in Theorem~\ref{thm4.3} can be dropped. Examples of infinite families of arithmetic reflection groups up to dimension $19$ have been given by Allcock \cite{Allc06}. These examples are obtained using the idea of \emph{doubling}: if a fundamental polyhedron $\hP$ of a reflection group $\Gamma$ has a face $\hF$ whose all dihedral angles with the other faces are equal to $\frac{\pi}2$, then we can double $\hP$ along $\hF$ to obtain a new polyhedron such that the reflections in its faces generate an index $2$ subgroup of $\Gamma$ (see Section~\ref{sec:examples} for the precise description of Allcock's redoubling procedure and related discussion). It is clear that the groups which are obtained by this procedure and its variations are non-maximal reflection groups. The necessity of the arithmeticity assumption is also well known. For example, consider the groups generated by reflections of the hyperbolic plane in the sides of a hyperbolic triangle with angles $\frac{\pi}2$, $\frac{\pi}3$, $\frac{\pi}m$ ($m\ge 7$). These groups are known to be maximal discrete subgroups of $\mathrm{Isom}(\Hy^{2})$ but all except finitely many of them are non-arithmetic. A similar construction is available for  hyperbolic $3$-space with triangular prisms replacing triangles (see \cite[Section~10.4.3]{MaclReid03}).

\medskip

Finally, let us mention that the spectral method that we described in this section has other interesting variations and applications. We refer to Peter Sarnak's lecture notes \cite{Sarnak14} for an enlightening discussion of some related topics.

\section{Effective results obtained by the spectral method} \label{sec:effective}

The proofs from the previous section can be made effective meaning that we can essentially enumerate all the possible candidates for the reflection groups in Theorems~\ref{thm4.1} and \ref{thm4.2}. In order to do so we need to look at the quantitative side of the finiteness theorems of Borel and Wang. The key ingredient for the quantitative analysis is provided by the results on minimal volume arithmetic hyperbolic $n$-orbifolds \cite{Bel04, Bel07, BelEme} (see also \cite{Bel14icm}) and the methods from these papers. In particular, all these results were obtained using the volume formula of G.~Prasad \cite{Pra89} which will also play a prominent role in our discussion.

More precisely, the quantitative analogue of Theorem \ref{thm4.2} is:
\begin{prop} (cf. \cite[Prop.~4.1]{Bel11})\label{prop5.1}
There are no cocompact congruence reflection subgroups in $\isom$ for $n \ge 13$, and no congruence reflection subgroups in $\isom$ for $n \ge 28$.
\end{prop}
It is proved by substituting in \eqref{eq:pr3} the precise lower bounds for the minimal volume. With some case-by-case considerations it should be not hard to bring the second bound down to $22$, so that it would agree with the Esselmann's result \cite{Esselmann96} bounding the dimension of non-cocompact arithmetic reflection groups.

A finer analysis based on the volume computations leads to the quantitative version of Theorem~\ref{thm4.2}:
\begin{prop} (cf. \cite[Prop.~4.2]{Bel11})\label{prop5.2}
The degrees of the fields of definition of cocompact congruence reflection subgroups of $\PO(n,1)$ are bounded by $6$ and their discriminants
satisfy the conditions in Table~\ref{tab52}.
\end{prop}
\begin{table}[!ht]
\caption{\label{tab52} The bounds for $\D_k$ depending on the dimension $n$ and the degree of the field $d = [k:\Q]$.}
$
\begin{array}{r|ccccc}
 & d = 2 & 3 & 4 & 5 & 6  \\
\hline
\rule{0pt}{2.5ex}
n = 4 & \le\!262 & \le\!2\,244 & \le\!19\,210 & \le\!164\,442 & \le\!1\,407\,650 \\
5 &  \le\!214 & \le\!1\,928 & \le\!17\,302 & \le\!155\,272 & \le\!1\,393\,406 \\
6 &  \le\!28 & 49,\ 81 \\
7 &  \le\!39 & \le\!205 & \le\!1\,062 \\
8,9 & \le\!13\\
10,11 & 5,\ 8\\
12 & 5 \\
\end{array}
$
\end{table}

Similar method can be employed to look at the other invariants (the Hasse--Witt symbol and the determinant of the quadratic form $f$ defining $\Gamma$) with an objective to give a list containing all the congruence-reflective quadratic forms. However, this has not been done yet. In small dimensions the list will be very large but for higher $n$ its size will reduce quickly --- compare with the possible fields of definition in Proposition~\ref{prop5.2}. Producing the list of the quadratic forms is a feasible task which would provide important data for the potential classification of arithmetic hyperbolic reflection groups. We shall come back to this discussion later in Section~\ref{sec:problems}.

The main issue about the results of the propositions is that they depend on an extra assumption --- the maximal reflection groups have to be congruence. We do not know yet how much, if any, information we loose by imposing this condition in dimensions $n\ge 4$ but for $n = 2$ and $3$ we do have the Lakeland's examples showing that not all arithmetic maximal reflection groups have this property \cite{LakelandPhD, Lakeland12}. Fortunately, there is a way to prove effective results in small dimensions without restricting to the congruence subgroups. Let us now review this method.

The idea goes back to \cite{LMR06} and \cite{Agol06icm}. From the point of view based on the conformal volume we can explain it as the \emph{Agol's trick:}

Let $\Gamma_0 < \isom$ be a maximal arithmetic subgroup containing a maximal reflection group $\Gamma$. Since a conjugation of a reflection in $\isom$ is again a reflection, the subgroup $\Gamma$ is normal in $\Gamma_0$. By the Vinberg's lemma (cf. Lemma~\ref{lem:max_refl}) we have that the quotient group $\Gamma_0/\Gamma$ is isomorphic to a finite subgroup $\Theta < \Gamma_0$ which is the group of symmetries of the Coxeter polyhedron $\hP$ of $\Gamma$. Consider $\Hy^n \subset \Sph^n$ embedded conformally as the upper half space of $\Sph^n$, so that $\isom$ acts conformally on $\Sph^n$. Normalize so that $\Theta$ acts isometrically on $\Sph^n$. Clearly, $V_{PC}(m, \OO) = V_{PC}(m, \hP/\Theta)$, where $\OO = \Hy^n/\Gamma_0$ and so $|\OO| = |\hP/\Theta|$. The orbifold embedding $\hP/\Theta \subset \Sph^n/\Theta$ is a conformal embedding, hence by one of the basic properties of the conformal volume we have $V_{PC}(n, \hP/\Theta) \le V_{PC}(n, \Sph^n/\Theta)$. The key observation is that we can give a good upper bound for $V_{PC}(n, \Sph^n/\Theta)$ if we manage to embed $\Theta$ in a finite reflection group $\Theta'$. Indeed, in this case we have
$$ V_{PC}(n, \Sph^n/\Theta) \le [\Theta':\Theta]\cdot V_{PC}(n, \Sph^n/\Theta') \quad \text{and} \quad V_{PC}(n, \Sph^n/\Theta') = \vol(\S^n),$$
which leads to
\begin{equation}\label{eq:eq4}
\vol_c(\OO) \le [\Theta':\Theta]\cdot \vol(\S^n).
\end{equation}
The required embedding $\Theta \hookrightarrow \Theta'$ is easy to obtain for $n = 2$ with the index $[\Theta':\Theta] \le 2$, Agol checked in his paper that for finite subgroups of $\mathrm{O}(3)$ we have $[\Theta':\Theta] \le 4$ and this resolves the case $n = 3$, it is also possible to extend this result to a more general class of quasi-reflective groups in dimension $n = 3$ which allowed their classification in \cite{BelMcl13}. In all other cases the classification of finite subgroups of $\mathrm{O}(n)$ is either not known or much more involved and we do not know how to bound the conformal volume of their quotients. Hence so far we can apply this trick for bounding the conformal volume only in dimensions two and three.

We now can substitute \eqref{eq:eq4} in \eqref{eq:Li-Yau} and use the bounds for the minimal volume and the spectral gap of the congruence quotients --- recall that according to Lemma~\ref{lem:max_ar} the maximal arithmetic subgroups are always congruence.

In the last section of his paper Agol indicated the possible quantitative implications of the method but he was missing some non-trivial technical ingredients required to make it work. It was observed later in \cite{Bel09} that one can combine Agol's method with the important technical results of Chinburg and Friedman \cite{ChinbFried86} in order to obtain the quantitative bounds. It was shown there, in particular, that the degree of the field of definition of arithmetic reflection groups in dimension $3$ is bounded above by $35$. In a joint work with Linowitz \cite{BelLin14}, we improved this bound to $9$, which essentially allows to give a list of all possible fields of definition (see \cite{BelLin14} for the details). The case of $n = 2$ was considered by Maclachlan in \cite{Macl11}. Summarizing the results we have:
\begin{theorem} \cite{Macl11, BelLin14}
The fields of definition of arithmetic hyperbolic reflection groups in dimension $2$ have degree at most $11$, and in dimension $3$ at most $9$.
\end{theorem}
This theorem complements Proposition~\ref{prop5.2} in a stronger form, as it does not impose any additional congruence hypothesis. The cited papers also provide explicit upper bounds for the discriminants of the fields of definition.

\section{Classification results in small dimensions} \label{sec:classification}

Consider a \emph{quadratic space} $(V,f)$, where $V$ is a finite dimensional vector space over a totally real number field $k$ and $f$ is a non-degenerate quadratic form on $V$. An $\cO$-module $L$ is called a \emph{quadratic lattice} if $L$ is a full rank $\cO$-lattice in $(V,f)$. A quadratic lattice is called \emph{even} if for the inner product associated with $f$ we have $(v,v)\in 2\cO$ for all $v\in L$, and \emph{odd} otherwise. The \emph{dual} $L^*$ of a quadratic lattice $L$ is the set of all vectors in $V$ having integer inner product with all vectors in $L$. A lattice is called \emph{unimodular} if $L = L^*$, in general, if the inner product is integral on $L$, we have $L \subseteq L^*$ and $\Delta(L) = L^*/L$ is a finite abelian group. It is called the \emph{discriminant group} of $L$, and its order is the \emph{determinant} $\mathrm{det}(L)$. A quadratic lattice $L$ is called \emph{strongly square-free} if the cardinality of the smallest generating set of $\Delta(L)$ as an $\cO$-module is at most $\frac12\mathrm{rank}(L)$ and every invariant factor of $\Delta(L)$ is square-free. If $k = \Q$ and the inner product associated with $f$ is $\Z$-valued on $L$, we shall call $L$ \emph{integral}. The \emph{level} of an integral lattice $L$ is defined to be the minimal positive integer $N$ such that $N(v,v)/2\in \Z$ for all $v\in L^*$. Integral quadratic lattices of signature $(n,1)$ or $(1,n)$ (and rank $n+1$) are called \emph{Lorentzian}. We refer to \cite{Conway-Sloane} for more material about quadratic lattices and related structures.

The group $\Gamma_L = \aut(L)$ of the automorphisms of an integral Lorentzian lattice is by the definition an arithmetic subgroup of the orthogonal group $\Ort(n,1)$. It can be shown that the automorphism group of a non strongly square-free lattice is always contained in the automorphism group of a strongly square-free one (cf. \cite{Allc12}).

Lorentzian lattices and their groups of automorphisms arise naturally in $K3$ surface theory, structure theory of hyperbolic Kac--Moody algebras, and many other fields. The question of their reflectivity was studied by Vinberg, Nikulin, Scharlau, Allcock, and others. There are also some related investigations about reflectivity of lattices defined over quadratic fields. In this section we shall review the classification results which come from this study. Except for an important work of Nikulin on $2$-reflective lattices discussed at the end, the other papers deal only with the lattices of small rank.

\subsection*{\textit{Nikulin, Allcock, Mark}} The case $n = 2$ is the first towards the general classification program. Although this case is easier than higher dimensions, the spectral method indicates that it is here that we can expect to encounter the largest number of examples of reflective lattices. In an important paper published in 2000 \cite{Nikulin00}, Nikulin classified the rank $3$ reflective strongly square-free Lorentzian lattices. He obtained a list of $1097$ lattices which fall into $160$ duality classes (a $p$-dual of an integral lattice $L$ is the sublattice of $L^*$ corresponding to the $p$-power part of $\Delta(L)$, it can be seen that $L$ and its $p$-dual have the same automorphism group). Note that since every lattice canonically determines a strongly square-free lattice, this classification does contain all integral Lorentzian lattices whose reflection groups are maximal under inclusion. The project was continued more recently by Allcock, who classified \emph{all} reflective integral Lorentzian lattices of signature $(2,1)$ \cite{Allc12}. He showed, in particular, that there are $8595$ such lattices which correspond to $374$ different reflection groups that fall into $39$ commensurability classes. He also checked that the $1097$ strongly square-free lattices previously obtained by Nikulin are contained in his list. The method is based on an analysis of the shape of the Coxeter polygons, which allows to reduce the list of candidates for the reflective lattices to a manageable size, and a subsequent application of Vinberg's algorithm.

It is worth mentioning that the TeX source file for the Allcock's paper (available from the mathematics e-prints arxiv, see \texttt{arXiv:1111.1264}) can be run as a Perl script that prints out all $8595$ lattices with their Gram matrices and other related data in computer-readable format. This brings the results to the form that is suitable for potential applications.

In her PhD thesis \cite{Mark_thesis}, Mark studied the classification of rank $3$ reflective lattices over quadratic extensions of $\Q$. To this end she adapted the method of Allcock to the real quadratic setting and used a different algorithm for checking reflectivity of a quadratic form termed the \emph{walking algorithm}. Mark's modifications to Allcock's method are inspired by the previous work of Bugaenko \cite{Bug84, Bug90, Bug92}.
The key difference between the walking and Vinberg's algorithms is the search space: whereas Vinberg's algorithm searches for the new vectors inside an $n$-dimensional polygonal cone, walking has a much more restricted searching area. The main result of \cite{Mark_thesis} is a complete classification of the rank 3 strongly square-free reflective arithmetic hyperbolic lattices defined over $\Z[\sqrt{2}]$. Mark showed that there are $432$ such lattices and provided their detailed description including the structure of the reflection groups. The methods developed in \cite{Mark_thesis} can be applied to the classification problem over other fields, which would be a natural next step of the classification project.

\subsection*{\textit{Scharlau--Walhorn.}}
In \cite{SchWal92}, the authors gave two explicit lists of maximal non-cocom\-pact arithmetic reflection groups in dimensions $n = 3$ and $4$. The groups are defined by the reflective integral quadratic lattices that are strongly square-free and isotropic. The lists in \cite{SchWal92} contain $49$ and $42$ lattices, respectively. Later Walhorn found that one example was missing from the list for $n = 4$, so there are in total $43$ such lattices \cite{Walhorn_thesis}. In the notation of \cite{SchWal92}, the 43rd lattice has the shape

\begin{quote}
$\Hy \perp \langle 1,7,7 \rangle$, where $\Hy$ is the even unimodular lattice of signature $(1,1)$ and $\langle a,b,c \rangle$ denotes the lattice of the diagonal quadratic form $ax^2+by^2+cz^2$. It has the determinant $D = -49$ and $r = 48$ fundamental roots.
\end{quote}

The lattices are shown to be reflective by Vinberg's algorithm, and as a bi-product of the algorithm application the authors also obtained various geometric invariants of the corresponding Coxeter polyhedra, such as the number of faces, the number of cusps, and other (see also \cite{Sch_prep, Walhorn_thesis} for more data). The papers do not provide Coxeter diagrams although they could have been produced from the algorithm output.

The authors indicate how to prove the completeness of the lists (we remark that for $n = 3$ they restrict to the isotropic quadratic forms and hence obtain only non-cocompact arithmetic reflection groups defined over $\Q$); the details of the proof for $n = 4$ are given in the dissertation of the second author \cite{Walhorn_thesis}.

\subsection*{\textit{Belolipetsky--Mcleod.}}
The previous enumeration for $n = 3$ is closely related to the study of reflective Bianchi groups.  For a square-free positive integer $m$ denote by $O_m$ the ring of integers of the imaginary quadratic field $\Q(\sqrt{-m})$. The \emph{Bianchi group} $Bi(m)$ is defined by $Bi(m)=\PGL(2, O_m)\rtimes \langle\tau\rangle$, where $\tau$ acts on $\PGL(2,O_m)$ as complex conjugation. The groups $Bi(m)$ can be regarded in a natural way as discrete subgroups of the group of isometries of the hyperbolic $3$-space $\Hy^3$. They are non-cocompact arithmetic subgroups of $\operatorname{Isom}(\mathbb{H}^3)$. One can also define the maximal discrete extension of $Bi(m)$ in $\operatorname{Isom}(\mathbb{H}^3)$, which is called the \emph{extended Bianchi group} and denoted $\widehat{Bi}(m)$. Reflectiveness of the Bianchi groups and their extensions was studied by many authors starting from the classical paper of Bianchi \cite{bianchi1891}. A triple of papers by Vinberg, Shaiheev, and Shvartsman published all in the same volume \cite{yaroslavl1988} made an important contribution to the topic by approaching it from the more general perspective of the Vinberg's program \cite{Vinb90, shaiheev90, shvartsman90}. In a related paper \cite{ruzmanov90}, Ruzmanoz considered an extended notion of reflectivity called \emph{quasi-reflectiveness} and gave the first examples of quasi-reflective Bianchi groups. We shall discuss this extension more carefully in Section~\ref{sec:quasirefl}. The research on these topics was concluded in \cite{BelMcl13}, where we proved the following classification result:

\begin{theorem}\label{thm:bianchi}
We have
\begin{enumerate}[(i)]
\item The Bianchi group $Bi(m)$ is reflective if and only if $m \le 19$, $m \neq 14, 17$.
\item The extended Bianchi group $\widehat{Bi}(m)$ is reflective if and only if $m \le 21$, $m = 30$, $33$, $39$.
\item The Bianchi group $Bi(m)$ is quasi-reflective if and only if $m=14$, $17$, $23$, $31$ and $39$.
\item The only quasi-reflective extended Bianchi groups are $\widehat{Bi}(23)$ and $\widehat{Bi}(31)$.
\end{enumerate}
\end{theorem}

In the proof, the finite list of candidates was produced using the spectral method and the final step of detecting the reflection groups was again performed by means of Vinberg's algorithm.
The paper also provides the Coxeter diagrams and other data for the reflection subgroups. Let us note that the list of the reflection groups in \cite{BelMcl13} is contained in the Scharlau--Walhorn classification for $n = 3$ but does not coincide with it because some integral quadratic forms give rise to non-cocompact arithmetic subgroups commensurable but not contained in Bianchi groups. It is
not hard to identify precisely the reflective groups from Theorem~\ref{thm:bianchi} in the table for $n=3$ in \cite{SchWal92}.

In his paper \cite{shaiheev90}, Shaiheev drawn the schematic pictures of the fundamental polyhedra of the reflection groups that he obtained (we note that there are some small mistakes in \cite{shaiheev90} and refer to \cite{BelMcl13} for the corrections). Now, with a complete classification available, it would be good to have a set of computer generated images of the Coxeter polyhedra of these reflection groups. Some nice examples of this type of polyhedra in the upper-half space model of $\Hy^3$ are presented in \cite{JJKSS15}. Another possible approach is to extend to the non-compact finite volume polyhedra the computer implementation of Andreev's theorem developed by Roeder \cite{Roeder07}.

\subsection*{\textit{Scharlau--Blaschke, Esselmann, Turkalj.}}
In all the above classification results the crucial step of determining the reflection subgroup is carried out by means of Vinberg's algorithm. Another approach to classification of the reflective integral lattices is based on the following lemma, which is also due to Vinberg:
\begin{lemma} (cf. \cite[Lemma 1.3]{SchWal92}) \label{lemma:genus}
Consider an integral lattice $L = \Hy \perp M$, where $\Hy$ is the even unimodular lattice of signature $(1,1)$ and $M$ is positive definite of rank at least $2$. If $L$ is reflective then the genus $\mathcal{G}(M)$ (which depends only on $L$) is totally reflective in the sense that every lattices $M' \in \mathcal{G}(M)$ is reflective.
\end{lemma}
We recall that in the positive definite case a lattice $M$ is called reflective if the reflection subgroup of its automorphisms group has no nonzero fixed vectors in $M\otimes\R$.

It is well known that under quite general conditions an integral lattice $L$ of signature $(n,1)$ does admit a decomposition of the form $\Hy \perp M$. In particular, this holds when $n \ge 4$ and $L$ is strongly square-free. Therefore, Lemma~\ref{lemma:genus} can be applied for the classification of non-cocompact reflection groups and allows to reduce the problem to the much better understood positive definite case. Esselmann proved in \cite{Esselmann96} that $20$ is the largest dimension for which there exist totally reflective genera. In \cite{SchBla96}, Scharlau and Blaschke used gluing technique to classify positive definite integral reflective lattices in dimensions $\leq 6$. Gluing theory mentioned here provides a method to construct general integral lattices that contain as a sublattice a direct sum of integral lattices of smaller dimension (see \cite[Chapter~4.3]{Conway-Sloane} for the details). In a recent preprint of Turkalj \cite{Turkalj15}, the work of Scharlau and Blaschke is combined with the other results to give a complete classification of the totally reflective primitive genera in dimension $3$ and $4$, which correspond to the Lorentzian lattices for the hyperbolic dimensions $n = 4$ and $5$, respectively. An explicit classification of the square-free totally positive genera for $n=4$ that appeared before in the Walhorn's dissertation \cite{Walhorn_thesis} is reproduced by Turkalj as a subset of the $n = 4$ case. The list in \cite{Turkalj15} contains $1234$ genera, of which $289$ are square-free and $52$ strongly square-free, in dimension $3$; and $930$ genera, of which $230$ are square-free and $88$ strongly square-free in dimension $4$. As it is expected, the total number of the reflective lattices decreases for bigger dimension.

\subsection*{\textit{Nikulin, Vinberg.}}
The case of $2$-reflective integral lattices is of a special interest because of its close connection with the theory of $K3$-surfaces. In particular, the classification of such lattices allows one to describe all algebraic complex surfaces of type $K3$ whose group of automorphisms is finite. Such a classification is now available thanks to the work of Nikulin and Vinberg.

Recall that an integral lattice $L$ is called \emph{$2$-reflective} if the subgroup of its group of automorphisms generated by $2$-reflections, i.e. the reflections whose primitive vectors have square $2$, is of finite index. Classification results for the $2$-reflective lattices go back to the first papers of Nikulin on the subject --- see \cite{Nikulin79}, \cite{Nikulin81a} and \cite{Nikulin84}. These papers cover all the cases except when the rank $r(L)$ is equal to $4$. The classification for the latter case was published by Vinberg only in 2007, although he obtained it as early as 1981 \cite{Vinb07}. The methods that are used in Nikulin's papers add algebraic geometry of $K3$-surfaces to the set of tools that we encountered before.

A complete list of $2$-reflective lattices can be found in the above cited papers, here we only reproduce the statistics --- see Table~\ref{tab6}. The $26$ lattices of rank $3$ are obtained in \cite{Nikulin84} (note that the lattices $S'_{6,1,2}$ and $S_{6,1,1}$ in the list there are isomorphic), the $14$ rank $4$ lattices are in \cite{Vinb07}, and the higher dimensions are treated in \cite{Nikulin81a}. For $r \ge 20$ the $2$-reflective integral lattices do not exist.

\begin{table}[!ht]
\caption{\label{tab6} $2$-reflective lattices.}
$
\begin{array}{c|ccccccccccccc}
r(L) & 3 & 4 & 5 & 6 & 7 & 8 & 9 & 10 & 11 & 12 & 13 & 14 & 15, \ldots, 19 \\
\hline
\rule{0pt}{2.5ex}
\#\text{ of lattices} & 26 & 14 & 9 & 10 & 9 & 12 & 10 & 9 & 4 & 4 & 3 & 3 & 1 \\
\end{array}
$
\end{table}

%
%

\section{Examples}\label{sec:examples}

\subsection*{\textit{Vinberg, Kaplinskaja}}
Hyperbolic reflection groups in dimensions $n \le 19$ were found by Vinberg \cite{Vinb72} and Vinberg--Kaplinskaja \cite{VK78}. They considered reflection subgroups of the groups of integral automorphisms of the quadratic forms
$$f(x_0, x_1, \dots, x_n) = -x_0^2 + x_1^2 + \ldots + x_n^2.$$
For $n = 2$ the form was investigated by Lagrange, Gauss, and later by Fricke. In his paper,  Fricke showed that the form $-x_0^2+x_1^2+x_2^2$ is reflective and described its Coxeter triangle fundamental domain \cite[pp. 64--68]{Fricke91}\footnote{I thank John Ratcliffe for suggesting this reference.}. The case $n =3$ first appeared, among other things, in the paper by Coxeter and Whitrow \cite{CoxWhit50}. For $n \le 17$ the form was investigated in \cite{Vinb72}, and the remaining $n = 18$ and $19$ were considered in \cite{VK78}. The Coxeter polyhedron in dimension $19$ is the most complicated one, it has $50$ faces and its symmetry group is isomorphic to the symmetric group $S_5$. More details can be read from the Coxeter diagrams that are presented in \cite{Vinb72} and \cite{VK78} for each of the cases (see also \cite[Chaper~28]{Conway-Sloane}).
In \cite{Vinb75}, Vinberg showed that the form $f$ is not reflective for $n \ge 25$ and indicated that the same should hold for $n\ge 20$ (see also \cite{VK78}). The proof of reflectivity in each of the cases is obtained by means of Vinberg's algorithm while non-reflectivity is shown by detecting infinite order elements in the quotient $\Gamma/\Gamma_r$ of the automorphism group by its reflection subgroup.

In \cite{Vinb72}, Vinberg also investigated the reflection groups of the quadratic forms
$$f_2(x_0, x_1, \dots, x_n) = -2x_0^2 + x_1^2 + \ldots + x_n^2.$$
He found that the form is reflective for $n \le 14$, non-reflectiveness of $f_2$ for bigger $n$ is confirmed in Mcleod's thesis \cite[Section~3.1.4]{mcleod_thesis}. The Coxeter diagrams for the reflection subgroups in the reflective case are given in Vinberg's paper.

\subsection*{\textit{Bugaenko}} By the Godement's compactness criterion, for $n \ge 4$ arithmetic groups defined by quadratic forms over $\Q$ are all non-cocompact. Thus in order to see cocompact higher dimensional examples we have to consider quadratic forms defined over the fields of degree at least $2$. This was first done by Bugaenko in 1980s and his examples still remain essentially the only ones of this type.

In \cite{Bug84}, Bugaenko investigated the reflection groups of the quadratic forms
$$f_{\sqrt{5}}(x_0, x_1, \dots, x_n) = -\frac{1+\sqrt{5}}{2}x_0^2 + x_1^2 + \ldots + x_n^2.$$
He proved that the form is reflective if and only if $n \le 7$. For $n = 8$ he found another admissible quadratic form over the field $\Q(\sqrt5)$ with discriminant $-(1+\sqrt{5})$ which is reflective \cite{Bug92}. This is the highest dimension for which we know examples of cocompact hyperbolic reflection groups.
The Coxeter diagrams for the reflection polyhedra are presented in \cite{Bug84} and \cite{Bug92}.

In another article \cite{Bug90}, Bugaenko considered the quadratic from
$$f_{\sqrt{2}}(x_0, x_1, \dots, x_n) =  -(1+\sqrt{2})x_0^2 + x_1^2 + \ldots + x_n^2,$$
which are shown to be reflective if and only if $n\le 6$. The Coxeter diagrams for $n \le 5$ are given in the paper (note that there are some minor typographical errors in the node-labeling for $n = 4$ and $5$, the vectors computed by the algorithm are correct), the data related to $n = 6$ provided by Bugaenko can be found in \cite[Tabels 2.1, 2.2]{Allc06}.

Reflectivity of the quadratic forms in each of the cases is checked by a variant of Vinberg's algorithm, with the modifications that make it work over the algebraic integers. To show the non-reflectivity Bugaenko systematically used a criterion of detecting an infinite order symmetry of the Coxeter polyhedron associated to a loxodromic isometry of $\Hy^n$.

Several other examples of cocompact arithmetic hyperbolic reflection groups similar to the ones that were considered by Bugaenko were found in Mcleod's thesis \cite{mcleod_thesis}. The case $n = 2$ over $\Q(\sqrt2)$ was thoroughly studied by Mark (see the previous section). In \cite[Section~4]{Allc13}, Allcock obtained an example of a cocompact reflection group in $\Hy^7$ from the reflection centralizer in Bugaenko's example in $\Hy^8$. This trick can be repeated to get an even more complicated example in $\Hy^6$. It would be interesting to know whether or not the resulting groups are commensurable with Bugaenko's examples.

\subsection*{\textit{Allcock, Potyagailo--Vinberg}}
In \cite{Allc06}, Allcock proved that there exist infinitely many finite-covolume (resp. cocompact) arithmetic hyperbolic reflection groups acting on hyperbolic space $\Hy^n$ for every $n \le 19$ (resp. $n \le 6$). This implies, in particular, that the maximality assumption in the finiteness Theorem~\ref{thm4.3} cannot be dropped. The construction is based on examples of Vinberg, Vinberg--Kaplinskaja and Bugaenko described above and a simple \emph{redoubling trick:}

Call a wall of a Coxeter polyhedron $\hP$ a doubling wall if the angles it makes with the walls it meets are all even submultiples of $\pi$. By the double of $\hP$ across one of its walls we mean the union of $\hP$ and its image under reflection across the wall. A polyhedron is called redoublable if it is a Coxeter polyhedron with two doubling walls that do not meet each other in $\Hy^n$. It is easy to show that the double of a Coxeter polyhedron $\hP$ across a doubling wall is itself a Coxeter polyhedron. Moreover, if the doubling wall is disjoint from another doubling wall, so that $\hP$ is redoublable, then the double is also redoublable. This allows one to iterate the procedure and produce an infinite series of finite volume Coxeter polyhedra.

The simplest redoublable polyhedra are the right-angled polyhedra, they have all the dihedral angles equal to $\pi/2$. These polyhedra were studied by Potyagailo and Vinberg in \cite{PVinb05}, who showed that they may exist only for $n \le 4$ in the compact case and for $n \le 14$ in the general finite volume case. Using a similar method, the second bound was later improved by Dufour to $n\le 12$ \cite{Dufour10}. Examples of compact right-angled polyhedra in $\Hy^n$ are known for all $n \le 4$ and finite-volume ones only for $n \le 8$ (see \cite{PVinb05}). For the other dimensions Allcock showed that many of the arithmetic examples discussed above are redoublable. It is worth mentioning that there exist also non-redoublable Coxeter polyhedra, the simplest example shown to me by Daniel Allcock is the hyperbolic triangle with all angles equal $\frac{\pi}{7}$ --- it is easy to check that the group generated by reflections in its sides \emph{does not have any} finite index reflection subgroups
except itself. It would be interesting to check if the same phenomenon occurs for the Borcherds polyhedron in $\Hy^{21}$ which is discussed below.

Allcock mentions that his method resembles Ruzmanov's construction of non-arithmetic Coxeter polyhedra from \cite{ruzmanov89}. The latter was recently further elaborated by Vinberg to produce some new examples of non-arithmetic hyperbolic reflection groups \cite{Vinb14}.

\subsection*{\textit{Mcleod, Mark}} The spectral method indicates that we should look for examples of arithmetic reflection subgroups in arithmetic lattices of small covolume. Recall that for every dimension the covolume of lattices in $\isom$ is uniformly bounded from below (by the Kazhdan-Margulis theorem \cite{KazMar68}), and the precise minimal value in the arithmetic case is known \cite{Bel04, Bel07, BelEme}. For most of the dimensions the minimum is attained on the arithmetic subgroups associated to the quadratic form $f$ considered by Vinberg, but quite surprisingly, for $n = 4k-1 \ge 7$, it corresponds to
$$f_3(x_0, x_1, \dots, x_n) = -3x_0^2 + x_1^2 + \ldots + x_n^2.$$
Reflection groups of these quadratic forms were investigated by Mcleod in \cite{mcleod11}. He showed that $f_3$ is reflective for $n \le 13$ and non-reflective for bigger $n$. The proofs use Vinberg's algorithm and some results of Bugaenko. Similar to the previous cases, arguably the most interesting example appears in the highest dimension $n = 13$, its Coxeter polyhedron has $22$ faces and the group of symmetries isomorphic to $\Z_2\times\Z_2$. The Coxeter diagrams are given in the Mcleod's paper.

The next natural step in this direction is to investigate the quadratic forms
$$f_m(x_0, x_1, \dots, x_n) = -mx_0^2 + x_1^2 + \ldots + x_n^2.$$
This was done by Mark for the case $m = p$, a prime number \cite{Mark15}. She showed that:
\begin{itemize}
\item[] $f_5$ is reflective for $2 \le n \le 8$;
\item[] $f_7$ and $f_{17}$ are reflective for $n = 2$ and $3$;
\item[] $f_{11}$ is reflective for $n = 2$, $3$, and $4$.
\end{itemize}
She also proved that for other $p$ and in higher dimensions $f_p$ is non-reflective. Together with some Nikulin's results for $n = 2$ this gives a complete list of the reflective forms of this type. Related results were also obtained by Mcleod in his thesis \cite{mcleod_thesis}, in particular, he gave a complete list of the reflective quadratic forms $f_m$ for all natural $m$ in all dimensions (see \cite[Table 3.1, page 37]{mcleod_thesis}).

\subsection*{\textit{Borcherds}} In \cite{Bor87}, Borcherds found an example of a non-cocompact arithmetic reflection group in $\Hy^{21}$, which was later shown by Esselmann \cite{Esselmann96} to have the largest possible dimension. Borcherds started from Conway's description of the group of automorphisms of the even unimodular Lorentzian lattice of rank $26$ in \cite{Conway83}: it is a semidirect product of the reflection subgroup (of infinite index) and the group of affine automorphisms of the famous \emph{Leech lattice}. We shall come back to this group in the next section. The finite volume Coxeter polyhedron $\hP^{21}$ in $\Hy^{21}$ comes out as a face corresponding to the spherical diagram of type $\Dn_4$ of the infinite volume $25$-dimensional Conway's polyhedron. A general method of determining the shape of a face of a Coxeter polyhedron motivated by this and other examples is described in \cite{Allc06}.

The polyhedron $\hP^{21}$ can also be obtained using Vinberg's algorithm. Following Bor\-cherds, its group is the reflection subgroup of $\Ort_0(g, \Z)$, where $g$ is a quadratic form of signature $(21,1)$ associated to the even sublattice $L$ of $\Z^{21,1}$ (i.e. $L$ consists of the integral vectors with the even sum of the coordinates). We can take the controlling vector $u_0 = v_0$, the first basis vector (assuming $(v_0, v_0) = -1$). Its stabilizer is generated by reflections corresponding to the remaining $21$ basis vectors, it is a finite Coxeter group of type $\Dn_{21}$. The next vector produced by the algorithm is $e_{22} = v_0 + v_1 + v_2 + v_3$, etc. The polyhedron $\hP^{21}$ has $210$ sides with $42$ of them corresponding to the $2$-reflections, and the remaining $168$ to the $4$-refections in $\Ort_0(g, \Z)$. It has a very large symmetry group isomorphic to $\PSL(3,\mathbb{F}_4) \cdot D_{6}$ of order $241920$ (here $D_{6}$ denotes the dihedral group of order $12$). It would interesting to try to draw its Coxeter diagram in a maximally symmetric way.

Daniel Allcock showed me a nice way to view the diagram of $\hP^{21}$ on the projective plane $\mathbb{F}_4\mathrm{P}^2$ over the field $\mathbb{F}_4$ with $4$ elements: one can index the faces of $\hP^{21}$ by the $21$ points, $21$ lines, and $168$ hyperconics in $\mathbb{F}_4\mathrm{P}^2$. The edges of the diagram are determined by the incidence relations. The resulting structure is invariant under the full group of automorphisms of $\mathbb{F}_4\mathrm{P}^2$ (including Galois conjugation and the point-line interchange), altogether producing the full group of symmetries of $\hP^{21}$. A cute algebraic geometric application of the Borcherds group related to this viewpoint was found by Dolgachev and Kondo \cite{DolgKond03}. They constructed a unique supersingular $K3$ surface in characteristic $2$ satisfying a set of equivalent properties whose automorphism group is the symmetry group of an infinite treelike polyhedron obtained by gluing together the copies of $\hP^{21}$.

The $21$-dimensional polyhedron considered above is a very special and quite complicated object, but the most complicated currently known finite volume example lives a few dimensions below. It was also discovered by Borcherds but in a different paper \cite{Bor00} and using a very different method. The idea of \cite{Bor00} is that many interesting reflection groups (in particular, most of the known examples in dimensions at least $5$) can be obtained from reflective singularities at cusps of modular forms of $\SL(2,\Z)$. This way Borcherds found new examples of arithmetic reflection groups without a priori writing down any roots and reflections! We shall review this method in the next section. The most complicated new example in $\Hy^n$ is described in \cite[p.~346]{Bor00}: it is a $17$-dimensional non-compact finite volume polyhedron with $960$ sides. Very little is currently known about geometry of this polyhedron.

\subsection*{\textit{Other examples}}
We conclude the discussion of classification and examples by mentioning brief\-ly some other results. There is another natural approach to the classification problem for hyperbolic reflection groups that rather than looking at the admissible quadratic forms and lattices in $\isom$ begins with analyzing the possible shapes of the Coxeter polyhedra in $\Hy^n$. The first class of the polyhedra that comes out here consists of the hyperbolic Coxeter simplices. Their study goes back to the work of Coxeter and Lann\'er in the first half of the XXth century. In 1950, Lann\'er enumerated bounded hyperbolic Coxeter simplices and showed that they exist only in dimensions $n \le 4$ \cite{Lan50}, later the enumeration was extended to the unbounded Coxeter simpleces of finite volume that exist in dimensions $n\le 9$. More recently, Johnson, Kellerhals, Ratcliffe and Tschantz described the commensurability classes of the hyperbolic Coxeter simplex reflection groups in all the  dimensions $9 \ge n\ge 3$ \cite{JKRT02}. They also showed that for $n \ge 4$ all these groups except for one $5$-dimensional example are arithmetic. In a series of papers Felikson and Tumarkin studied other types of the hyperbolic Coxeter polyhedra (without connection to arithmeticity). We refer to \cite{FelTum14} and the references therein for the related results.

\section{Reflective modular forms} \label{sec:refl forms}

Let $\Gamma$ be a lattice in $\SL(2, \R)$. A \emph{modular form} of weight $k$ with respect to $\Gamma$ is a complex-valued function $f$ on the hyperbolic plane $\Hy = \{z\in \C \mid \mathrm{Im}(z) > 0 \}$ in the upper half-plane model which is holomorphic on $\Hy$, holomorphic at the cusps of $\Gamma$ and satisfies the equation
\begin{equation*}
f\left(\frac{az+b}{cz+d}\right) = (cz+d)^k f(z), \text{ for all } z\in\Hy \text{ and } \begin{pmatrix} a & b \\ c & d \end{pmatrix} \in \Gamma.
\end{equation*}
A modular form is called a \emph{cusp form} if it vanishes at the cusps of $\Gamma$. For example, if $\Gamma = \SL(2, \Z)$ this condition means that $f(z) \to 0$ when $z \to i\infty$.

The theta (or Howe) correspondence assigns to a cusp form $f$ a cuspidal automorphic form $\phi_f$ on the bounded symmetric domain associated with a group $\Ort(m,n)$. We are interested in the singular theta correspondence, which allows $f$ to have poles at the cusps and as the output produces a meromorphic modular form $\phi_f$. It can be shown then that under certain conditions on $f$ the singularities of $\phi_f$ are reflection hyperplanes of an arithmetic reflection group or a quasi-reflection group. We shall proceed with a more precise description of the correspondence.

Let $L$ be an integral lattice of signature $(m,n)$ and $\Gr(L)$ denote the Grassmannian of the maximal ($m$-dimensional) positive definite subspaces of $L\otimes\R$. It is a symmetric space of dimension $mn$ acted by the orthogonal group $\Ort(m,n)$. We note that in this section the Lorentzian lattices have signature $(1,n)$ which gives opposite signs of inner products compared to the rest of the paper. We decided not to change the notation in order to comply with the literature. Given an element $v \in \Gr(L)$ and $\lambda \in L\otimes\R$, we denote by $\lambda_{v^+}$ and $\lambda_{v^-}$ the projections of $\lambda$ onto the positive definite space represented by $v$ and the negative definite space orthogonal to $v$.

Suppose the lattice $L$ is even. The Siegel theta function of a coset $L+\gamma$ of $L$ in $L^*$ is
$$
\Theta_{L+\gamma}(\tau; v) = \sum_{\lambda \in L+\gamma} \exp(\pi i\tau\lambda^2_{v^+} + \pi i\bar{\tau}\lambda^2_{v^-}),
$$
where $\tau \in \Hy$ and $v \in \Gr(L)$. Combining these for all elements of $L^*/L$ gives a $\C[L^*/L]$-valued function called the \emph{Siegel theta function} of $L$:
$$
\Theta_{L}(\tau; v) = \sum_{v \in L/L^*} e_\gamma \Theta_{L+v}(\tau; v),
$$
where $e_\gamma$ denote the elements of the standard basis of the group ring $\C[L^*/L]$.

When dealing with theta functions of lattices, half-integral weight vector-valued modular forms naturally occur. They can be defined using the metaplectic double cover $\widetilde{\SL}(2, \R)$. We refer to \cite{Bor98a} for the details of this construction. The Siegel theta function $\Theta_{L}(\tau; v)$ is a vector-valued modular form of weight $(m/2,n/2)$ and type $\rho_L$, where $\rho_L$ is the Weil representation of the group $\widetilde{\SL}(2, \Z)$ on the vector space $\C[L^*/L]$.

Let $F$ be another vector-valued modular form which has weight $(-m/2,-n/2)$ and type $\rho_L$. Then the product $F(\tau)\overline{\Theta(\tau;v)}$ is a modular form of weight $0$. If this product is of sufficiently rapid decay at $i\infty$ (which occurs if $F$ is a cusp form), we can take the integral
\begin{equation*}
\Phi_F(v) = \int_\mathcal{F} F(v)\overline{\Theta(\tau;v)} \frac{dx dy}{y^2},
\end{equation*}
where $\mathcal{F} = \{\tau \in \Hy \mid \ |\tau| \ge 1,\ |\mathrm{Re}(\tau)|\le 1/2| \}$ is the usual fundamental domain for $\SL(2, \Z)$. This gives us a function $\Phi_F$ on $\Gr(L)$ invariant under a congruence subgroup of $\aut(L)$. The map $F(\tau) \to \Theta_F(v)$ is essentially the original theta correspondence.

Suppose now that we allow $F(\tau)$ to have singularities at the cusps but require it to be holomorphic on $\Hy$. The integral above diverges wildly. Harvey and Moore used ideas from quantum field theory to show that it is still possible to make sense of the integral by the \emph{regularization} \cite{HM96}. Their construction was further generalized by Borcherds in \cite{Bor98a}. The idea of the regularization is to truncate the integration domain in such a way that most of the wildly non-convergent terms vanish. The remaining non-convergent terms are of polynomial growth and can be dealt with easily. The truncated domains are
$$
\mathcal{F}_t = \{\tau \in \Hy \mid \ |\tau| \ge 1,\ |\mathrm{Re}(\tau)|\le 1/2|,\ \mathrm{Im}(\tau) \le t \}.
$$
The regularized value of the integral is defined as the value at $r = 0$ of the analytic continuation of
\begin{equation*}
\lim_{t\to\infty} \int_{\mathcal{F}_t} F(\tau)\overline{\Theta(\tau;v)} \frac{dx dy}{y^{2+r}}.
\end{equation*}
This defines a more general map $F(\tau) \to \Phi_F(v)$ which is called the \emph{singular theta correspondence}. It is easy to check that the singularities of $\Phi_F(v)$ occur on sub-Grassmannians of the form $v^{\perp}$ for $v\in L^*$, $(v.v)<0$, where there is a nonzero coefficient corresponding to $v$ in the Fourier expansion of $F$ at the cusp. If the singularities occur along the reflection hyperplanes of the underlying lattice $L$, we shall call the modular form $F(\tau)$ a \emph{reflective modular form}. Not all reflective lattices correspond to reflective modular forms but many particularly interesting examples do have this property.

The construction can be generalized to modular forms of the groups $\Gamma$ different from $\SL(2, \Z)$, moreover, it is often possible to use scalar valued modular forms of level $N$ instead of the vector valued modular forms. We refer to \cite{Bor00} for the details. A useful sufficient condition for a modular form to be reflective is given in \cite[Lemma~11.2]{Bor00}. For example, if $N$ is a square-free integer and $\Gamma = \Gamma_0(N) = \left\{ \begin{pmatrix} a & b \\ c & d \end{pmatrix} \in \SL(2, \Z) \mid\ c\equiv 0 \text{ mod } N\right\}$ is a congruence subgroup, a modular form $F(\tau)$ for $\Gamma$ of weight $\frac{m-n}2$ is reflective for an even lattice $L$ of signature $(m,n)$ and level $N$ if the poles of $F(\tau)$ at all cusps of $\Gamma$ are simple.

\medskip

We conclude this section with some examples of reflective modular forms from \cite{Bor00}.

\begin{example} \label{ex:8.1}
The first case to consider is $N = 1$, $\Gamma = \SL(2, \Z)$ and $L$ is an even integral lattice of signature $(m,n)$ and level $1$. It is well known that the modular forms of $\Gamma$ form a polynomial ring generated by the Eisenstein series $E_4(\tau) = 1 + 240q + 2160q^2 + \ldots$ and $E_6(\tau) = 1- 504q - 16632q^2 - \ldots$, with $q = e^{2\pi i \tau}$, of weights $4$ and $6$, respectively. The dimensions of the spaces of modular forms of different weights are given by the coefficients of the Hilbert function $1/(1-x^4)(1-x^6)$. We refer to \cite{Miyake} for these and other related facts from the classical theory of modular forms. The first weight in which we have a non-trivial modular form is $k = 12$, and the critical form is $\Delta(\tau) = \eta(\tau)^{24} = q \prod_{k>0}(1 - q^k)^{24}$. The forms $f = E_4(\tau)^k/\Delta(\tau)$ have simple poles at the cusp of $\Gamma$ at $i \infty$. It follows that they are reflective modular forms for even lattices $L$ of level $N = 1$ and signature $m-n \ge -24$, $m-n = 0 \text{ mod }8$.

Here are some concrete cases:

For $m - n = -24$ we can take $f = 1/\Delta(\tau) = q^{-1} + 24 + 324q + \ldots$ of weight $-12$ and a simple pole at the cusp. Two examples arising here are of a special interest: in the Lorentzian case we have a lattice $L = II_{1, 25}$, which is a quasi-reflective lattice first discovered by Conway \cite{Conway83}; in signature $(2,26)$ the reflective lattice $II_{2, 26}$ plays an important role in the arithmetic mirror symmetry studied by Gritsenko and Nikulin.

For $m - n = -16$ we take $f = E_4(\tau)/\Delta(\tau) = q^{-1} + 264 + 8244q + \ldots$ of weight $-8$ and find the reflective Lorentzian lattice $II_{1, 17}$. Similarly, for $m - n = -8$ with $f = E_4(\tau)^2/\Delta(\tau)$ $= q^{-1} + 504 +73764q + \ldots$  of weight $-4$ we obtain the reflective lattice $II_{1, 9}$. The arithmetic reflection groups associated with these lattices were first described by Vinberg in \cite{Vinb75}.
\end{example}

\begin{example} \label{ex:8.2}
Suppose the level $N = 2$. The group $\Gamma = \Gamma_0(N)$ has two cusps which can be taken as $i \infty$ and $0$. The ring of modular forms for $\Gamma$ is a polynomial ring on generators $-E_2(\tau) + 2E_2(2\tau) = 1 + 24q + 24q^2 + \ldots$ of weight $2$ and $E_4(\tau)$ of weight $4$. The Hilbert function is $1/(1-x^2)(1-x^4)$.

As $N$ is square-free, all poles of order at most $1$ are reflective by \cite[Lemma~11.2]{Bor00}. There are also other possible reflective singularities but we will not consider them here. By looking at the form $\Delta_{2+}(\tau)^{-1} = \eta(\tau)^{-8}\eta(2\tau)^{-8}$ of weight $-8$ with simple poles at the cusps, we see that all level $2$ even lattices of signature at least $-16$ have reflective modular forms. The Lorentzian lattices $II_{1,17}(2^{+8})$ and $II_{1,17}(2^{+10})$ are quasi-reflective as was the case for the lattice $II_{1, 25}$ in the previous example (we refer to \cite{Conway-Sloane} for the notation used here). The next example $L = II_{1,17}(2^{+6})$ gives us the $17$-dimensional arithmetic reflection group discovered by Borcherds in \cite{Bor00}. This example was mentioned at the end of the previous section.
\end{example}

\section{More about the structure of the reflective quotient} \label{sec:quasirefl}

Let $\Gamma_0 < \isom$ an arithmetic subgroup, $\Gamma \mathrel{\unlhd} \Gamma_0$ its maximal subgroup generated by reflections in hyperplanes, and $\Theta = \Gamma_0/\Gamma$ the reflective quotient. We have the following possibilities for the group $\Theta$:
\begin{itemize}
\item[(a)] finite group;
\item[(b)] affine crystallographic group;
\item[(c)] non-amenable group.
\end{itemize}

In case (a) the group $\Gamma$ is an arithmetic reflection group, while in (b) and (c) it has infinite covolume and hence is not a lattice. Case (b) is known as \emph{quasi-reflective} or \emph{parabolic-reflective}. The second term refers to the fact that in this case the group $\Theta$ is virtually isomorphic to an affine group generated by parabolic transformations of $\Hy^n$.  Recall that a discrete group is called \emph{amenable} if it has a finitely-additive left-invariant probability measure. Finite groups and affine crystallographic groups are amenable hence
the first two cases can be joined together into the \emph{amenable type}. This type is relatively rare, the known results imply that generically the reflective quotients are non-amenable.

Note that by the Tits alternative these are the only possible cases. Indeed, the group $\Theta$ is a finitely generated linear group hence by \cite{Tits72} it is either virtually solvable or contains a non-abelian free subgroup. A virtually solvable group acting discretely on hyperbolic space (recall that by Lemma~\ref{lem:max_refl}, the group $\Theta$ is isomorphic to the group of symmetries of the Coxeter polyhedron of $\Gamma$) is virtually free abelian, which brings us to a one of the first two cases. From the other hand, it is well known that the groups that contain non-abelian free subgroups are non-amenable. It is worth pointing out that the case (b) can happen only if $\Gamma_0$ is not cocompact, in cocompact case we have a simple alternative that either the reflective quotient $\Theta$ is finite or it is non-amenable.

\medskip

Most of this survey is dedicated to the reflective groups (a). The first and arguably the most interesting example of a quasi-reflective group was discovered by Conway in \cite{Conway83}. He showed that for the group of automorphisms of the $26$-dimensional even unimodular Lorentzian lattice its reflective quotient is isomorphic to the group of affine automorphisms of the Leech lattice. Therefore, we have an example of a quasi-reflective arithmetic group in hyperbolic dimension $n = 25$. The proof in Conway's paper is very short but it relies on a lot of results from the previous study of the Leech lattice by Conway and others. It is conjectured that $n = 25$ is the largest dimension where there exists a quasi-reflective group. This conjecture can be possibly resolved by Esselmann's method \cite{Esselmann96}, which he applied to find the maximal dimension of isotropic reflective lattices, but such a proof is not available so far. In his doctoral dissertation \cite{Barnard_thesis}, Barnard showed that the conjecture can be deduced from an open conjecture of Burger, Li and Sarnak about the automorphic spectra of orthogonal groups \cite{BurgSarn91, BurgLiSarn92}\footnote{I thank Richard Borcherds for mentioning to me this important result.}. It is interesting to note that another conjecture considered in \cite{BurgSarn91} has already appeared in this survey while we were discussing the spectral method, Barnard's approach is based on reflective modular forms and highlights yet another relation between arithmetic reflection groups and the spectrum of the Laplacian.

In \cite{Nikulin96}, Nikulin proved using his method that in any dimension $n$ there exists only finitely many maximal arithmetic hyperbolic quasi-reflective groups (see Theorem~1.1.3 ibid.). He mentions that it is not difficult to show that there are no any such groups for $n \ge 43$, and also conjectures that the sharp bound should be given by the Conway's example. It would be interesting to find a spectral proof for the finiteness theorem, which may potentially give better quantitative bounds for these groups in a fixed dimension. The only result of this kind available so far can be found in \cite[Section~5]{BelMcl13}, where we used the Agol's trick (cf. Section~\ref{sec:effective}) and the classification of the plane crystallographic groups to obtain a good quantitative bound for the maximal quasi-reflective groups in dimension $n = 3$. Note that similar to the reflective case, the maximality assumption is essential: an example of an infinite sequence of arithmetic quasi-reflective groups in dimension $2$ is given in \cite[Example 1.3.4]{Nikulin96}.

Following Conway, quasi-reflective groups related to the Leech lattice were studied by Borcherds who, in particular, found several other examples in smaller dimensions (see \cite[Theorem~3.3]{Bor90}). Later examples of quasi-reflective groups were constructed using quasi-reflective modular forms \cite{Bor00, Barnard_thesis}. In \cite{ruzmanov90}, Ruzmanov considered quasi-reflective groups in dimension $3$ from the geometric viewpoint. He introduced the notion of a quasi-bounded Coxeter polyhedron and found examples of quasi-reflective Bianchi groups. This research was concluded in \cite{BelMcl13}, where all the quasi-reflective Bianchi groups and extended Bianchi groups are classified (cf. Theorem \ref{thm:bianchi}(iii, iv)).

\medskip

By the classification of the possible types of the reflective quotients, the upper bounds for the dimension of the reflective and quasi-reflective groups imply the lower bound for the dimension in which the quotient is necessarily non-amenable. A different approach to the problem was undertaken by Meiri in \cite{Meiri14}, who directly constructed non-abelian free subgroups of the reflective quotients in sufficiently large dimension for the forms defined over $\Q$. A part of Meiri's argument is closely related to the proof of non-reflectivity of the high dimensional quadratic forms in \cite{Vinb84}. Unfortunately, the quantitative bounds on the dimension obtained in \cite{Meiri14} are far from sharp.

\medskip

It would be interesting to know what we can say about the reflective quotient group when it is non-amenable: Is it (relatively) hyperbolic? $\mathrm{CAT}(0)$? Has uniformly bounded spectral gap?.. etc.

\section{Open problems}\label{sec:problems}

\subsection*{\textit{Generalizations}}
In the very beginning of the paper we discussed Petrunin's question \cite{MOver}, we recall it again here:

\begin{problem}\label{prob9.1}
Do there exist any hyperbolic lattices in the spaces of large dimension which are generated by elements of finite order?
\end{problem}
As before, the question can be restricted to the arithmetic hyperbolic lattices. The answer is unknown in both cases, but we can expect it to be \emph{negative}.

A related problem appears in a recent paper by Fuchs, Meiri and Sarnak (cf. \cite[page~1621]{FMS14}):
\begin{problem}\label{prob9.2}
Are there any hyperbolic lattices generated by reflections and Cartan involutions (also called ``reflections in points'') in the hyperbolic spaces of sufficiently large dimension?
\end{problem}
A negative answer to this question would allow to settle Conjecture~2 in \cite{FMS14}. An example of a lattice generated by Cartan involutions in $\Hy^8$ can be found in \cite[Theorem~5.3]{Allc99}.
This problem is, of course, a very special case of Problem~\ref{prob9.1}.

The groups of isometries of the hyperbolic $n$-space has real rank $1$.
It is worth mentioning that for the irreducible lattices in higher real rank semisimple Lie groups $H$ the situation is very different.  We can consider a lattice $\Gamma_0$ in such a group $H$ which has a non-central element $g$ of finite order. Let $\Gamma$ be the normal subgroup of $\Gamma_0$ generated by all the conjugates of $g$. By the Margulis normal subgroup theorem (see \cite[Chapter~IV]{Marg91}), the group $\Gamma$ is then itself a lattice in $H$. We can leave as an exercise for the reader to check that it is generated by a finite set of $g$-conjugates. Hence we can easily produce various examples of lattices in such groups $H$ which are generated by elements of finite order.

Basic examples of this kind of lattices in higher rank groups come from the orthogonal groups $\Ort(n, m)$, $n, m \ge 2$. Similar to the case of the hyperbolic $n$-space (corresponding to $m = 1$), we can consider here lattices generated by reflections. For instance, for the quadratic form
$$f_{n,m}(x_1, x_2, \dots, x_{n+m}) = -x_1^2 - \ldots - x_m^2 + x_{m+1}^2 + \ldots + x_{m+n}^2$$
the group generated by reflections will be an arithmetic lattice in $\Ort(n, m)$ according to the arithmeticity and normal subgroup theorems of Margulis \cite{Marg91}. These lattices, however, will never be Coxeter groups because as the lattices in a higher rank Lie group they have Kazhdan's property~(T) \cite{Kazdan67}, while on the other hand the infinite Coxeter groups are known not to have Kazhdan's property \cite{BozJan88}. It would be good to know more about the algebraic structure of these groups:
\begin{problem}\label{prob9.3}
Obtain examples of presentations for lattices in $\Ort(n, m)$, $n, m \ge 2$, generated by reflections or, more generally, for irreducible lattices in higher rank semisimple Lie groups generated by elements of finite order.
\end{problem}

Considering the infiniteness of the family of higher rank reflection groups, the natural question is which of them are really interesting. In \cite{Bor00}, Borcherds suggested that \emph{interesting reflective lattices} should be associated to reflective modular forms and gave examples of such lattices (cf. Section~\ref{sec:refl forms}). In a recent preprint \cite{Ma15}, Ma studied the basic class of $2$-reflective modular forms proving that there are only finitely many $2$-reflective lattices of signature $(2, n)$ with $n \ge 7$ and there are no such lattice when $n \ge 30$. Here a lattice $L$ is called \emph{$2$-reflective} if the subgroup of its group of automorphisms generated by $-2$-reflections is of finite index and $L$ is associated with a reflective modular form. Clearly, the second condition is crucial for the finiteness result. The largest $n$ for which we know an example of this kind is $n = 26$ \cite[p.~344]{Bor00} with the corresponding lattice $L = II_{2, 26}$ (cf. Example~\ref{ex:8.1}).

\subsection*{\textit{Towards classification}}
There are two main problems that appear on the way towards classification of arithmetic hyperbolic reflection groups:
\begin{problem}\label{prob9.4}
Find good bounds for the arithmetic invariants of the reflective quadratic forms in arbitrary dimension.
\end{problem}

\begin{problem}\label{prob9.5}
Check reflectivity of a given quadratic form.
\end{problem}

The quantitative bounds that can be extracted from the proofs of the finiteness theorems in \cite{ABSW08} or \cite{Nikulin07_fin} are huge and have no practical value. In Section~\ref{sec:effective} we explained how Problem~\ref{prob9.4} can be solved under a certain additional arithmetic assumption (requiring that the maximal reflection groups are congruence) or for dimensions $n \le 3$. One can try to push the technique from the low dimensions to higher $n$, or try to investigate what kind of limitations are in fact implied by the congruence assumption. We can also produce the conditional list and use it as a heuristic for generating all the compatible examples. It is plausible that this list would actually cover all higher dimensional examples. For sufficiently large $n$, say $n \ge 10$, the conditional list of the candidates would be quite small.
Another interesting direction is to try to find all lattices that are associated to reflective modular forms. Is there any connection between these lattices and the congruence reflection groups?

The main tool for checking reflectivity of a quadratic form is Vinberg's algorithm discussed above. There are several computer implementations of this algorithm but non of them is publicly available or standard. The case when the ring of integers $\cO$ of the defining field is not a PID requires a special attention as in this case the transformations from $\Gamma$ are not necessarily given by the matrices with $\cO$-entries. I do not know examples of arithmetic reflection groups that would highlight this issue but they may exist, in particular, in the spaces of small dimension.

\subsection*{\textit{About examples}}
Geometry of the high dimensional hyperbolic Coxeter polyhedra remains mysterious in many ways. This refers also to the known examples of such polyhedra. Here we can state the following problem:
\begin{problem}\label{prob9.6}
Investigate geometrical properties and combinatorial structure of the high dimensional hyperbolic Coxeter polyhedra.
\end{problem}
In this statement the high dimension may refer to $n \ge 6$ in the compact case and to $n \ge 10$ for the finite volume non-compact hyperbolic polyhedra. There are only handful known examples in these dimensions (cf. Section~\ref{sec:examples}). Which of them are redoublable or mixable in the sense of \cite{Allc06} and \cite{Vinb14}, respectively? What are the covering/comensurability relations between the higher dimensional examples? How many faces, vertices and cusps do these polyhedrons have? These numbers are known for most of the examples but are there any interesting relations between them beyond the ones that we already know?

\subsection*{\textit{A hyperbolic version of the Cartan--Killing classification}} My attention to this problem was brought by Daniel Allcock. A fundamentally important feature of the finite and affine Coxeter groups is their connection to Lie theory. This connection extends to hyperbolic Coxeter groups and Kac--Moody theory. In a series of papers Gritsenko and Nikulin isolated the key properties of Borcherds' fake monster Lie algebra and stated the corresponding classification problem for the class of the Generalized Kac--Moody algebras which they call the \emph{Lorentzian Kac--Moody algebras} (see \cite{GrNik02}).

The classification problem can be stated in terms of root systems. Let us call a set $\Pi$ of spacelike (positive-norm) vectors in $\Eucln1$ a \emph{simple root system} if $(v, v')$ is non-positive and lies in $\frac12(v,v)\Z$, for all $v, v' \in \Pi$. The integral span $L$ of the simple roots is called the \emph{root lattice}. The Weyl group $W$ means the group generated by the reflections in the roots.

\newpage

A simple root system $\Pi$ spanning $\Eucln1$ satisfies the \emph{Gritsenko--Nikulin conditions} if:
\begin{itemize}
 \item[(i)] there exists a Weyl vector $\rho \in \Eucln1$ such that $(v, \rho) = -(v,v)/2$ for all $v \in \Pi$;
 \item[(ii)] the normalizer of $W$ has a finite index in the orthogonal group $\Ort(L)$;
 \item[(iii)] the group $W$ has the \emph{arithmetic type} (cf. \cite[Section 1.4]{Nikulin96} or \cite[p. 326--327]{Allc15} for several equivalent definitions of this notion).
\end{itemize}
In the terminology of \cite{Nikulin96} such root systems are said to have the \emph{restricted arithmetic type}. Under these conditions, it turns out that the Weyl vector $\rho$ is unique and must be timelike or lightlike (i.e. $(\rho,\rho) < 0$ or $= 0$, respectively). In the timelike case, the (projectivized) Weyl chamber is a finite volume polyhedron in hyperbolic space, while in the lightlike case, the Weyl chamber has infinite volume and infinitely many sides.

\begin{problem}\label{prob9.7}
Classify the simple root systems of rank at least $3$ that satisfy the Gritsenko--Nikulin conditions and corresponding Lorentzian Kac--Moody algebras.
\end{problem}

This problem is closely related to the classification problem for the Lorentzian reflective and quasi-reflective lattices. The main difference is the existence of the Weyl vector $\rho$ which is related to the automorphic features of the Lorentzian algebras (recall, in particular, the Borcherds reflective modular form mentioned in Section~\ref{sec:examples}). The known results about hyperbolic reflection groups allowed Nikulin to prove finiteness theorems for the Lorentzian Kac--Moody algebras (see \cite{Nikulin96}). Nikulin and Gritsenko constructed families of examples of such algebras and gave a partial classification of them in rank $3$ (see \cite{GrNik02} and the references therein). A complete classification of the rank $3$ root systems with a timelike Weyl vector that satisfy the Gritsenko--Nikulin conditions was obtained recently by Allcock \cite{Allc15}. The special features of the Gritsenko--Nikulin root systems make the classification problem more accessible, while their relation to Lie theory makes it particularly interesting.

\subsection*{\textit{Complex reflection groups}} A complex reflection is an automorphism of a finite dimensio\-nal complex vector space that fixes a complex hyperplane. In contrast with the real case, complex reflections do not necessarily have order $2$. Finite linear groups generated by complex reflections were classified by Shephard and Todd \cite{ShepTodd54}. They are distinguished from all other finite linear groups by the property that their algebras of invariants are free. Realizing Hermitian symmetric spaces as bounded symmetric domains in a complex vector space allows us to define complex reflections in such spaces. It is well known that the only bounded symmetric domains which admit totally geodesic complex hypersurfaces are the complex hyperbolic spaces $\C\Hy^n$ (corresponding to the groups $\U(n,1)$) and the domains of type IV (of the groups $\Ort(n,2)$). It follows that only these spaces admit complex reflections.

The study of lattices in the complex hyperbolic space $\C\Hy^n$ goes back to the nineteenth century. In the 1883 paper \cite{Picard83}, Picard investigated the lattices $\SU(2,1; O_d)$, where $O_d$ is the ring of integers of the imaginary quadratic field $\Q(\sqrt{-d})$, acting on $\C\Hy^2$. These groups are called the \emph{Picard modular groups}. In a recent article \cite{PaupWill13}, Paupert and Will showed that for $d=1$, $2$, $3$, $7$, $11$ the Picard modular groups are up to a finite index generated by real reflections (i.e. antiholomorphic involutions that have a real totally geodesic plane fixed). It is not known which of the Picard groups are complex reflective. In an influential paper \cite{Mostow80}, Mostow constructed Dirichlet fundamental domains for certain groups generated by complex reflections in $\C\Hy^2$, both arithmetic and non-arithmetic. The only other complex hyperbolic space in which a non-arithmetic lattice is known is $\C\Hy^3$ \cite{DelMost86}. Higher dimensional examples of arithmetic complex hyperbolic reflection groups were constructed by Allcock in \cite{Allc00a, Allc00b}. His method is close to the constructions in the real hyperbolic spaces considered in this survey. In \cite{Allc00a}, Allcock obtained examples of lattices generated by complex reflections in $\C\Hy^5$, $\C\Hy^9$ and $\C\Hy^{13}$ by considering the automorphism groups of Lorentzian lattices over the Eisenstein integers $O_3$. In \cite{Allc00b}, he used a related construction to give several other examples of lattices, including examples in $\C\Hy^4$ and $\C\Hy^7$ which do not appear on the list of Deligne and Mostow \cite{DelMost86}. The Allcock group in a record high dimension $13$ is again related to the Leech lattice. Some interesting examples of complex hyperbolic reflection groups are considered in \cite{Deraux06}, \cite{Stover14}, and \cite{DPP15}. A survey of the known constructions of complex hyperbolic lattices by Parker can be found in  \cite{Parker09}.

Currently available results about complex hyperbolic reflection groups and lattices essentially fall into several types of the known constructions. There are no general finiteness theorems or classification attempts. Most of the questions that were discussed in this survey about real hyperbolic reflection groups can also be asked in the complex case, but here the results will turn into conjectures or open problems. It is not even clear, for instance, if we should expect that the dimension of arithmetic reflection complex hyperbolic lattices is bounded from above.

The case of domains $\mathbb{D}_n$ of type IV is different. Here the group $\Ort(n,2)$ has real rank $2$, so we can obtain many examples of lattices generated by complex reflections for any $n$ using the construction for higher rank groups described in the beginning of this section. The natural question that comes out again is how to narrow this class. Vinberg suggested to use as a criterion freeness of the algebra of automorphic forms on $\mathbb{D}_n$ invariant under $\Gamma$ thus generalizing the free polynomial invariant algebras of the finite complex reflection groups. 
There are some particularly interesting examples of lattices generated by complex reflections arising this way. Igusa's paper \cite{Igu62} provides such an example for $n = 3$. His results were largely extended by Vinberg in \cite{Vinb10}, whose construction provides this type of examples of arithmetic complex reflection groups in the domains of type IV for $n = 4$, $5$, $6$, $7$.

In view of Borcherds work \cite{Bor00} discussed in Section~\ref{sec:refl forms}, one can ask if arithmetic complex reflection groups can be related to modular forms in a similar way it is done for the real reflective lattices. There are some known examples of this kind (see e.g. \cite{Allc00a}) but the general picture remains mysterious. 

\subsection*{\bf Acknowledgement.} I am grateful to Ernest Borisovich Vinberg and Daniel Allcock for their comments and corrections, which helped to significantly improve the quality of this paper. I would like to thank Alice Mark for sending me the preliminary version of her thesis and Igor Dolgachev for sending me his lecture notes. Comments from  Richard Borcherds, Benjamin Linowitz, John Parker,  John Ratcliffe, and the referee on the preliminary version of the paper were very helpful. This work was partially supported by the CNPq and FAPERJ research grants.

\bibliographystyle{alpha}
\bibliography{mbel}

\end{document}